\documentclass{gtart2}
\usepackage{amssymb}
\def\<{\langle}
\def\>{\rangle}
\def\|{{\ |\ }}
\newcommand{\implies}{\Longrightarrow}
\newcommand{\inv}{^{-1}}              
\newcommand{\basl}{\backslash}

\newcommand{\C}{{\mathbb C}}
\newcommand{\N}{{\mathbb N}}    
\newcommand{\Z}{{\mathbb Z}}
\newcommand{\R}{{\mathbb R}}
\newcommand{\bH}{{\mathbb H^2}}
\newcommand{\dbH}{{\overline{\mathbb H}^2}}
\newcommand{\MCG}{{\mathcal{MCG}}}
\renewcommand{\phi}{\varphi}
\input epsf
\def\relabelbox{%
  \hbox\bgroup%
}%
\def\endrelabelbox{%
\egroup%
}%
\def\relabel #1#2 {%
  \special{ps:/a {} def}%
  \smash{\rlap{#2}}%
}%
\def\adjustrelabel <#1,#2> #3#4 {%
  \special{ps:/a {} def}%
  \smash{\rlap{\kern #1 \raise #2\hbox{#4}}}%
}%
\def\extralabel <#1,#2> #3 {\smash{\rlap{\kern #1 \raise #2\hbox{#3}}}}%

\newtheorem{lemma}{Lemma}[section]
\newtheorem{theorem}[lemma]{Theorem}
\newtheorem{prop}[lemma]{Proposition}
\newtheorem{cor}[lemma]{Corollary}
\newtheorem{definition}[lemma]{Definition}
\newtheorem{lem}[lemma]{Lemma}

\newtheorem{remarks}[lemma]{Remarks}

\begin{document}


\title{Orderings of mapping class groups after Thurston} 

\author{
Hamish Short\\Bert Wiest}

\address{CMI, Universit\'e de Provence\\39 Rue Joliot Curie\\13453 Marseille,
FRANCE\\
\smallskip\\\rm Email:\stdspace\tt hamish@cmi.univ-mrs.fr\\
bertw@cmi.univ-mrs.fr}

\begin{abstract}
We are concerned with mapping class groups of surfaces with
nonempty boundary. 
We present a very natural method, due to Thurston, of 
finding many different left orderings of such groups. The construction 
involves equipping the surface with a hyperbolic structure, embedding
the universal cover in the hyperbolic plane, 
and extending the action of the mapping class group on it
to its limit points on the circle at infinity. We classify all orderings 
of braid groups which arise in this way. Moreover, restricting to a certain
class of ``nonpathological'' orderings, we prove that there 
are only finitely many conjugacy classes of such orderings.
\end{abstract}

\keywords{mapping class group, braid group, orderable group}

\primaryclass{20F36, 20F60}\secondaryclass{57M07, 57M60}

\makeshorttitle


We shall be concerned with surfaces $S$ with nonempty boundary and a
finite set of punctures, and their mapping class groups $\MCG(S)$, 
ie the group of homeomorphisms 
$S\to S$ which map $\partial S$ identically and permute the punctures, 
up to isotopy. It was first proved by Dehornoy \cite{D} that
braid groups (ie mapping class groups of punctured disks) are
left-orderable. A topological proof of this result
was given in \cite{FGRRW}, and the extension to mapping class groups of 
general surfaces with boundary can be found in \cite{RW}. Here we present 
a very natural method, due to Thurston \cite{Thur}, of 
finding many different left orderings of such groups. 
In brief, one equips the
surface with a hyperbolic structure, lifts it to $\bH$, attaches to 
this cover its limit points on the circle at infinity, and notices that
there is a natural action of the mapping class group on the (circular)
boundary of the resulting space which fixes a point, and thus an action 
on $\R$. 
We classify the set of orderings of braid groups which arise from Thurston's
construction; more precisely, we divide these orderings into two disjoint 
classes, which we call orderings of finite, respectively infinite, type; 
the orderings inside each of the classes are classified by combinatorial
means. Finite type orderings are discrete, and they occur in only finitely
many conjugacy classes. By contrast, there are uncountably many infinite
type orderings, and all of them are non-discrete.

The outline of the paper is as follows. In the first section we give a 
short introduction to orderable groups and survey some known results about 
them. In the second section we present Thurston's construction. 
In the third section we define finite and infinite type 
orderings, and state our classification theorems. 
Section four to six are concerned with finite type orderings:
in section 4 we describe a different method of constructing
orderings, using ``curve diagrams''. In the section 5 we prove
that the set of orderings arising from curve diagrams 
is very easy to understand and classify. Moreover, we prove that up to
conjugacy only a finite number of orderings arise in this way.
In the sixth section we prove the classification theorems for finite
type orderings. 
The strategy is to associate to every point of $\R$ with orbit of finite type
a curve diagram such that the  orderings arising from this point and from 
the curve diagram agree. Thus we obtain, via curve diagram orderings, a 
good understanding of Thurston type orderings. In section 7 we prove
the results about the infinite type case.


\section{Orderable groups}

In this section we define orderable groups and survey some known results
about them.

\begin{definition} \sl
A group $G$ is {\it left orderable\/} (resp.\ {\it right orderable\/}) 
if there is a total order $<$ on $G$ which is  invariant under left 
mutiplication (resp.~ right multiplication),
that is, such that, for all $a,b \in G$, $a<b$, $a=b$ or $b<a$,
and for all $g\in G$, $a<b$ implies that $ga<gb$ (resp. $ag<bg$).

A group $G$ is {\it two--sided orderable\/}, or {\it bi-orderable} 
if there is a total order on $G$ which is respected by multiplication on 
the left and  multiplication on the right: 
i.e. $a<b\implies ga<gb$ and $ag<bg$.

Two left orderings $<$ and $\prec$ on a group $G$ are {\it conjugate}
if there exists a $g\in G$ such that $a\prec b$ if and only if
$ag < bg$. So two left orderings are conjugate if ``one is obtained from
the other by right translation in the group''.
\end{definition}

\begin{remarks}\rm
(1) The following observation will be extremely important in what follows.
If a group $G$ acts on the left by orientation preserving homeomorphisms 
on $\R$, then every point $x$ in $\R$ with free orbit (ie $Stab(x)=\{1_G\}$)
gives rise to a left ordering on $G$, by defining $g>h :\hspace{-1mm}
\iff g(x)>h(x)$.
We have for every $f\in G$ that $fg(x)>fh(x) \iff g(x)<h(x)$, since
the action of $f$ preserves the orientation of $\R$; this implies that the 
ordering is indeed left invariant. Note that different points in $\R$ may
give rise to different orderings.
If a point $x$ does not have free orbit, it still gives rise to a
{\it partial} left-invariant ordering.

(2) In fact, a countable group is left-orderable if and only if it has an 
action by orientation preserving homeomorphisms on $\R$ such that only the
trivial group element acts by the identity-homeomorphism, see for instance
\cite{GhysLima}. 

(3) A left orderable group is torsion--free: if an element $x$ had order $n$, 
and if $1<x$, then it would follow that $1<x<x^{2}<\dots < x^{n-1}<x^n =1$.


(4) The ``positive cone'' of the ordering, $P=\{g\in G \mid g>1\}$
has the properties that $G=P\sqcup\{1\}\sqcup P^{-1}$, and that 
$PP\subset P$. Conversely, given a subset with these two properties, 
a left order $<$ can be defined by $a<b :\hspace{-1mm}\iff a^{-1}b\in P$.
Similarly, a right order $\prec$ is obtained from 
$a\prec b :\hspace{-1mm}\iff ab^{-1}\in P$.
(In particular, a group is left orderable if and only if it is 
right-orderable.) The orders are total because of the frst property, 
and transitive because of the second.
The orders are bi-orders if and only if we have in addition that 
$g\inv P g \subseteq P$ for all $g\in G$.

(5) The following classes of groups are bi-orderable: (a) finitely generated 
torsion--free abelian groups;  (b) finitely generated free groups (this 
is a result of Magnus, see eg\ \cite{RK}); (c) more generally, residually
free groups, like fundamental groups of closed surfaces (this is due
to Baumslag, see \cite[p.484]{S}).

(6) The standard reference for orderings on groups is Rehmtulla and Mura's 
book \cite{RM}. 
\end{remarks}
 
We now give three examples of attractive results about orders on groups.

(1) In 1966, Neville Smythe \cite{S} used the orderability of
surface groups to prove
that  any null--homotopic curve on a surface $S$ is the image under 
projection of an embedded unknotted loop in $S\times I$.

(2) In 
\cite[Question N]{N1}, L. Neuwirth asked whether all knot groups 
are two sided orderable. In his review of the state of the questions 
several years later \cite{N2}, 
he thanks N. Smythe for pointing out that the trefoil
knot group (which is isomorphic to the braid group on three strings 
$B_{3}$), is not two sided orderable.
To show this, consider the presentation $B_{3}=\<x,y; x^{2}=y^{3}\>$.
In this group one has $xy\neq yx$.
Suppose that $B_{3}$ is two sided orderable, and without loss of 
generality suppose that $yx>xy$. 
Conjugating this inequality by $x$ 
yields
$$
xy=x(yx)x\inv > x(xy)x\inv = x^2yx\inv = y^4x\inv = yx^2x\inv = yx,
$$
so we have the contradiction that $xy>yx$.

Neuwirth reformulates the question as `Are knot groups 
left orderable?'.
A positive answer to this question follows from an observation
by J.~Howie and H.~Short \cite{HS} that knot groups are locally 
indicable (every non--trivial finitely generated subgroup has $\Z$
as a homomorphic image), together with a theorem of Burns and Hale 
\cite{BH} that locally indicable groups are left orderable.
The converse of Burns and Hale's theorem is known to be false
- see \cite{B} and \cite[Theorem 5.3]{FGRRW}.

(3) A conjecture of Kaplansky, called the Zero Divisor Conjecture, asserts
that if $R$ is a ring without zero divisors and $G$ is a torsion--free group 
then the group ring $RG$ has no zero divisors. The hypothesis that $G$ be
torsion free is necessary, for if $G$ contains an element $x$ of order $n$ 
then $(1-x)(1+x + \dots + x^{n-1})=0$ in $RG$. 
The conjecture is known to hold for left orderable groups. In fact, is 
not hard to
see that left orderable groups have the ``two unique product'' property
which implies that the conjecture holds for them (see eg \cite{Pa}, and also
Delzant \cite{Del} and Bowditch \cite{Bow} for some recent remarks 
about this property.)


\section{Orderings of mapping class groups using hyperbolic geometry}

In this section we present the definition of orders on mapping class
groups of surfaces which we learned from W Thurston, and prove that
they all extend the subword-ordering of Elrifai- Morton. 
The idea comes from the following classical situation. As is well known, 
every closed surface of genus $g\geqslant 2$ can carry a hyperbolic 
structure; ie there is a homeomorphism between the universal cover $S^\sim$ 
of $S$ and the hyperbolic plane $\bH$ such that the covering transformations 
are isometries of $\bH$. 
There is a natural closure $S^\eqsim\cong \dbH$ of $ S^\sim\cong \bH$, 
defined by adding the so-called {\it circle at infinity} $S^1_\infty=
\partial \dbH$. 
Points of this circle can be defined as classes of geodesics, or 
quasi-geodesics, $\gamma: [0,\infty) \to \bH$, staying a bounded distance 
apart. 
The covering action of $\pi_1(S)$ on $S^\sim$ extends to an action
on $S^\eqsim$. So in particular, we have an action of $\pi_1(S)$
on the circle at infinity by homeomorphisms; this action has been much 
studied (for a good modern exposition of all this see \cite{Ghys}). 
Even stronger, every homeomorphism of the surface lifts and 
extends to a homeomorphism of $S^\eqsim$; however, there is a 
$\pi_1(S)$-family of possible choices of lift, and therefore we get no 
well-defined action of $\MCG(S)$ on $S^1_\infty$.

Instead of closed surfaces, Thurston considers surfaces $S$ with nonempty 
boundary, a finite number of punctures, and $\chi(S)<0$. 
Again, one can obtain a hyperbolic structure on $S$ in which $\partial S$ is
a geodesic and the punctures are cusps; this is done by identifying 
$S^\sim$ with a subset of the hyperbolic plane, in such
a way that the covering translations are isometries of $\bH$. 
The boundary of this subset is just the union of the lifts of $\partial S$;
in particular it is a union of geodesics in $\bH$, and it follows that 
$S^\sim$ is convex in the hyperbolic metric. Moreover, the set of limit 
points of $S^\sim$ on the circle at infinity $\partial \dbH$ is a Cantor set
in $\partial \dbH$.
The closure $S^\eqsim$ of $S^\sim$ in $\dbH$, ie $S^\sim$ with its limit 
points on the circle at infinity attached, is homeomorphic to a closed disk;
$\partial S^\eqsim$ is a circle, also containing $S^\eqsim \cap \partial \dbH$ 
as a Cantor set. 

We now fix, once and
for all, a basepoint of $S^\sim$ anywhere on $\partial S^\sim$. 
We denote the component of $\partial S^\sim$ which contains the base
point by $\Pi$. The basepoint projects to a basepoint of $S$ in $\partial S$,
and $\Pi$ is an infinite cyclic cover of one component of $\partial S$. 
We consider the set of geodesics in $S^\eqsim$ starting at the basepoint - 
they are parametrized by the interval $(0,\pi)$, according
to their angle with $\Pi$. 
We shall denote by $\widetilde{\gamma}_\alpha$ the geodesic with angle 
$\alpha\in (0,\pi)$ and by $\gamma_\alpha$ its projection to $S$.
Since $S^\sim$ is hyperbolically convex, each 
point of $\partial S^\eqsim$ can be connected to the basepoint by
a unique geodesic (possibly of infinite length) in $S^\eqsim$, and for 
points in $S^\eqsim \basl \Pi$ this is one of the geodesics 
$\widetilde{\gamma}_\alpha$ with $\alpha \in (0,\pi)$.
This construction proves
\begin{lem} \sl There is a natural homeomorphism between 
$\partial S^\eqsim \basl \Pi$ and $(0,\pi)$. \end{lem}
As in the case of closed surfaces, we have an action of $\pi_1(S)$
on $S^\eqsim$, which restricts to an action on $\partial S^\eqsim$.
However, this time we have more:
\begin{prop} \sl There is a natural action by orientation preserving
homeomorphisms of $\MCG(S)$ on $\partial S^\eqsim \basl \Pi \cong (0,\pi)$.
\end{prop}
\begin{proof} Every homeomorphism $\varphi\co S \to S$ has a canonical
lift $\widetilde{\varphi}\co S^\sim \to S^\sim$, namely the one
that fixes the basepoint of $S^\sim$, and thus all of $\Pi$. Moreover, 
$\widetilde{\varphi}$ has an extension $\overline{\widetilde{\varphi}}\co 
S^\eqsim \to S^\eqsim$. The restriction of this homeomorphism to
$\partial S^\eqsim$ is invariant under isotopy of $\varphi$, and fixes $\Pi$,
and thus yields a well-defined homeomorphism of $\partial S^\eqsim \basl \Pi$.
\end{proof}
\begin{cor} \sl $\MCG(S)$ is left-orderable. \end{cor}
%

\newpage
\cl{
\relabelbox
\epsfbox{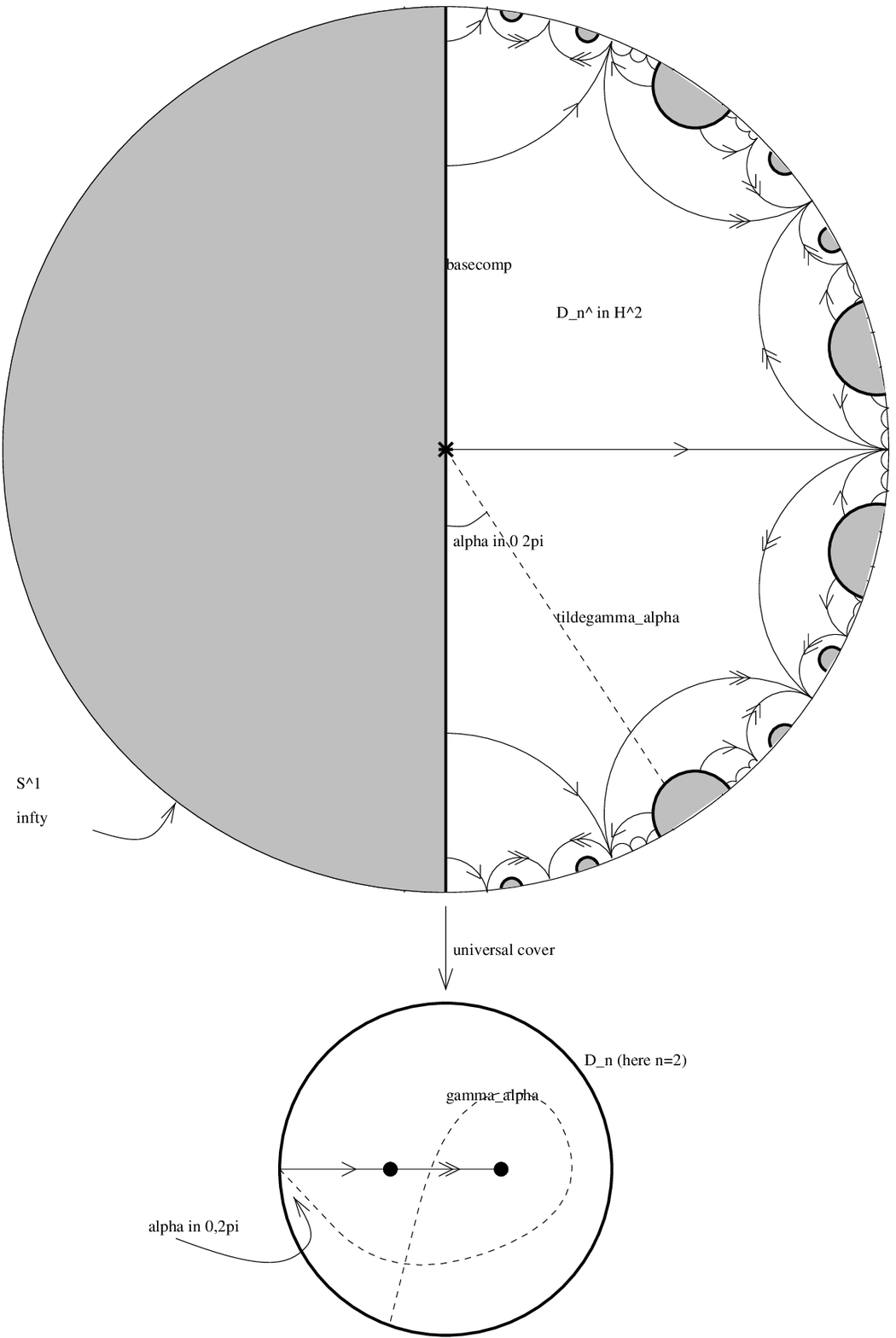} 
\relabel{D_n^ in H^2}{\small $D_n^\sim \subset \bH$}
\relabel{alpha in 0 2pi}{\small $\alpha \in (0,\pi)$}
\relabel{alpha in 0,2pi}{\small $\alpha \in (0,\pi)$}
\relabel{universal cover}{\small universal cover}
\relabel{basecomp}{\small $\Pi$}
\relabel{gamma_alpha}{\small $\gamma_\alpha$} 
\relabel{tildegamma_alpha}{\small $\widetilde{\gamma}_\alpha$} 
\relabel{D_n (here n=2)}{\small $D_n$ (here $n=2$)}
\relabel{S^1}{\small circle at} 
\relabel{infty}{\small infinity \thinspace$\partial\dbH$}
\endrelabelbox
}
\begin{figure}[htb]
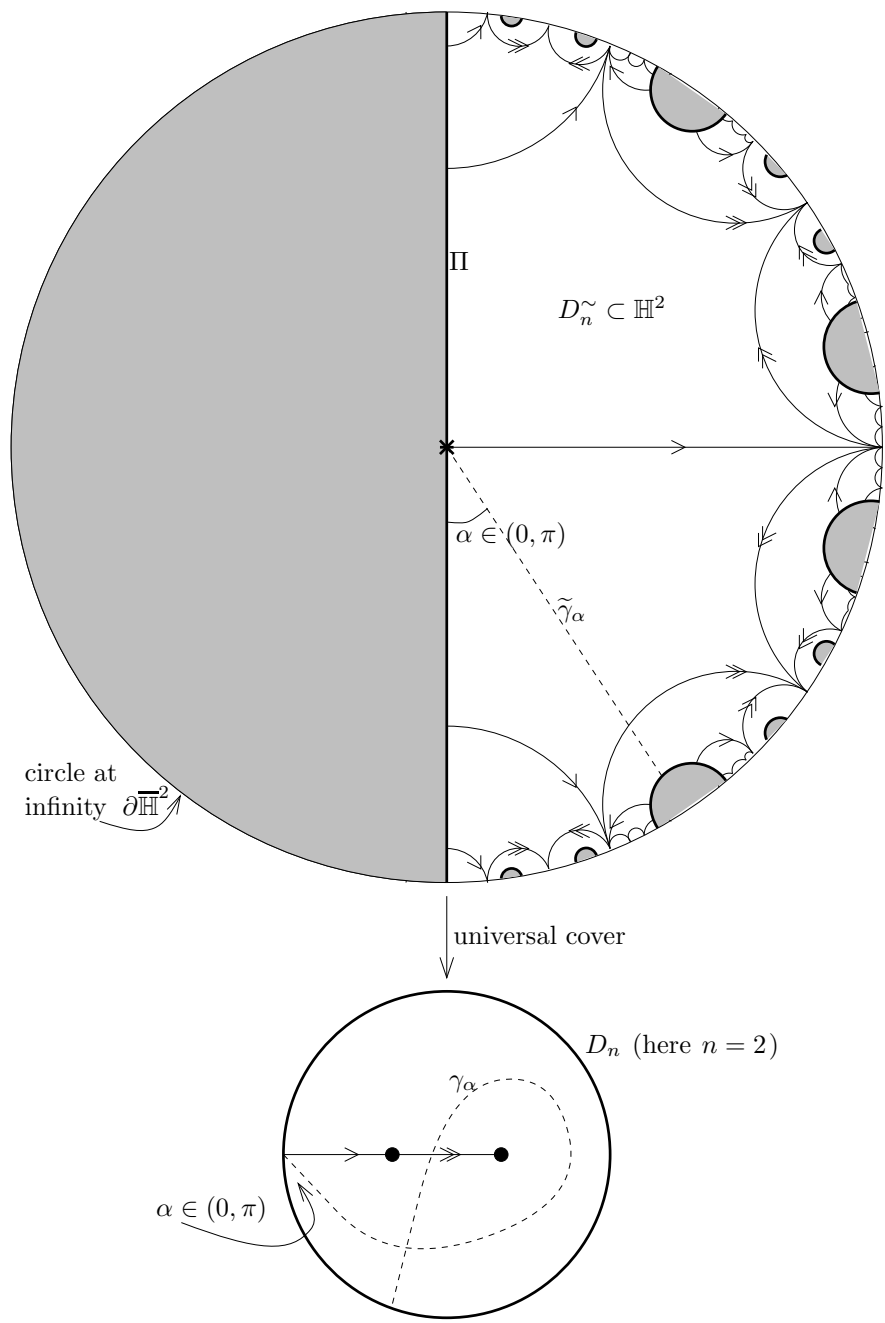

\caption{picture of $S^\sim$ in $\bH$ (here $S$ is a twice-punctured disk)}
\label{hyp}
\end{figure}
\newpage

%
\begin{proof} 
No nontrivial element of $\MCG(S)$ acts trivially on $(0,\pi)$, because 
if such an element existed, it would in particular fix
all liftings of the basepoint of $S$, and thus induce the 
identity-homorphism on $\pi_1(S)$, contradicting the hypothesis
that the element be non-isotopic to the identity.

The result now follows from Remark 1.2(2),
because $(0,\pi)$ is homeomorphic to $\R$. 
However, there is an elementary proof in our situation. We choose
arbitrarily a finite generating set of $\pi_1(S)$, and denote the
end points 
of the liftings of these elements by $s_1,\ldots, s_k \in (0,\pi)$.
A left order on $\MCG(S)$ is now defined inductively: if $\phi(s_1)>s_1$
then $\phi>1$ (and the same with $>$ replaced by $<$); if $\phi(s_1)=s_1$,
but $\phi(s_2)>s_2$, then $\phi>1$ as well, and so on; this is a total
order, because we have that $\phi(s_i)=s_i$ for all $i$ if and only if 
$\phi=1$.
\end{proof}

However, for the rest of the paper we shall be less interested in orderings
of this type, but rather in orderings induced by the orbits of single
geodesics, ie in orderings of the type introduced in Remark 1.2(1). 

\begin{prop} \sl For any simple closed geodesic $\tau$ in the surface $S$ 
we have for the positive Dehn-twist $T$ along $\tau$ that
$T(\alpha) \geqslant \alpha$ for any $\alpha \in (0, \pi)$. If
$\gamma_\alpha$ intersects 
$\tau$ at least once, then the inequality is strict.
\end{prop}
\begin{proof}
If $\gamma_\alpha$ is disjoint from $\tau$, then $T(\alpha)=\alpha$. 
If, on the contrary, $\gamma_\alpha$ 
intersects $\tau$, and hence any curve isotopic to $\tau$, any number of 
times (possibly infinitely often), then we denote by $T_i(\gamma_\alpha)$ 
($i\in \N$) the 
curve obtained from $\gamma_\alpha$ by applying the Dehn twist to the
first $i$ intersections of $\gamma_\alpha$ with $\tau$ and ignoring
all following intersections; we denote by $T_i(\alpha)$ its end point in 
$\partial D_n^\eqsim \basl \Pi$. We have 
$T(\alpha)=\lim_{i\to \infty} T_i(\alpha)$.

We now claim that $(T_i(\alpha))_{i \in \N}$ is a strictly increasing 
sequence. To simplify notation, we shall prove the special case 
$T_1(\alpha)>\alpha$, the proof in the general case is exactly the same. 
In the universal cover $D_n^\eqsim$ we consider the lifting of the curve 
$T_1(\gamma_\alpha)$: starting at the basepoint, it sets off along
$\widetilde{\gamma}_\alpha$, up to the first intersection with some
lifting $\widetilde{\tau}$ of $\tau$. There it turns right, walks along
$\widetilde{\tau}$ up to the next preimage of the intersection point,
where it encounters a different lifting ${\widetilde{\gamma}}'_\alpha$
of $\gamma_\alpha$. There it turns left, following this lifting all
the way to $\partial D_n^\eqsim \basl \Pi$. The crucial point now is
that $\widetilde{\gamma}_\alpha$ and ${\widetilde{\gamma}}'_\alpha$
intersect $\widetilde{\tau}$ at the same angle, because the two
intersections are just different liftings of the same intersection
between $\gamma_\alpha$ and $\tau$ in $D_n$. It follows that  
$\widetilde{\gamma}_\alpha$ and ${\widetilde{\gamma}}'_\alpha$ do not
intersect, not even at infinity, for if they did they would determine
a hyperbolic triangle in $D_n^\eqsim$ two of whose interior angles 
already add up to 180 degree, which is impossible. This implies the claim, 
and thus proves the proposition.
\end{proof} 
\begin{cor} \sl All total orderings of the braid group $B_n$ considered in 
this paper extend the subword ordering of Elrifai-Morton \cite{EM, W}.
More precisely, if a curve $\tau$ in $D_n$ encloses a precisely twice
punctured disk and $T^{1/2}$ is the half-Dehn-twist along $\tau$ 
interchanging the two punctures then $T\circ \phi > \phi$ for any
$\phi \in B_n$ and any ordering $>$ of Thurston-type.
\end{cor}

\begin{proof} It suffices to prove that
$T^{1/2}(\alpha)\geqslant \alpha$ for all $\alpha \in (0,\pi)$. 
If there existed an $\alpha\in (0,\pi)$ with $T^{1/2}(\alpha)<\alpha$
then it would follow that 
$T(\alpha)=T^{1/2}\circ T^{1/2}(\alpha) < T^{1/2}(\alpha) < \alpha$
(where the first inequality holds since $T^{1/2}$ is orientation
preserving), in contradiction with the proposition. \end{proof}


\section{Main results}

Let $D_n$ be the closed unit disk in $\C$, with $n$ punctures lined up in
the real interval $(-1,1)$ ($n\geqslant 2$).
We shall mainly be interested in the case $S=D_n$, where the mapping class 
group is a braid group $\MCG(D_n)=B_n$.
We recall that for $\alpha \in (0,\pi)$ we denote by $\gamma_\alpha$ the
geodesic which starts at the basepoint with angle $\alpha$ with 
$\partial D_n$, and by $\widetilde{\gamma}_\alpha$ its preimage in the 
universal cover starting at the basepoint of $D_n^\sim$.

\begin{definition}\sl 
A geodesic $\gamma_\alpha$, $\alpha \in (0,\pi)$, is said to be of 
{\it finite type} if it satisfies at least one of the following conditions

(a) there exists a 
finite initial segment $\gamma_\alpha^t$ such that any two punctures that
lie in the same path component of $S\basl \gamma_\alpha^t$ also lie in the
same path component of $S\basl \gamma_\alpha$, or

(b) it falls into a puncture, or

(c) it spirals towards a simple closed geodesic.

If a geodesic $\gamma_\alpha$ is not of finite type then we say it is of 
{\it infinite type}. We also define the ordering of $B_n$ induced by a 
geodesic $\gamma_\alpha$ to be of finite or infinite type if $\gamma_\alpha$
is of finite or infinite type.
\end{definition}

An infinite type geodesic looks as follows. All its self intersections occur
in some finite initial segment $\gamma_\alpha^t$. At least one of the path
components of $S\basl \gamma_\alpha^t$ contains three or more punctures in
its interior, and the closure of $\gamma_\alpha \basl \gamma_\alpha^t$
is a geodesic lamination without closed leaves inside such a component.
In particular, there is a pair of punctures which are seperated by the
whole geodesic, but not by any finite initial segment.

\begin{definition} For a geodesic $\gamma_\alpha$ of finite respectively 
infinite type we say that it {\it fills the surface in finite} respectively 
{\it infinite time} if
all punctures lie in different path components of $S \basl \gamma_\alpha$.
By contrast, a geodesic $\gamma_\alpha$ {\it does not fill} the surface if
$S\basl \gamma_\alpha$ has a path component that contains two punctures.
\end{definition}

The aim of the rest of the paper is to prove the following theorems.
Recall that every point $\alpha \in (0,\pi)$ gives rise to a - possibly
partial - ordering of $\MCG(S)$. The first theorem gives criteria for
these orderings to be total or, equivalently, for the orbit of $\alpha$
to be free.
\begin{theorem}\sl (a) If a geodesic $\gamma_\alpha$ does not fill $S$, 
then the orbit of $\alpha \in (0,\pi)$ is not free. 

(b) If $\gamma_\alpha$ is of finite type, then the converse holds as well: 
if $\gamma_\alpha$ fills the surface, then $\alpha$ has free orbit.

(c) The set of points $\alpha \in (0,\pi)$ such that $\gamma_\alpha$ 
is of infinite type consists of uncountably many, discrete, 
nonisolated points. All but countably many of those that fill $S$ in 
infinite time have free orbit.
\end{theorem}
The next theorem gives a classification of order types that
arise from Thurston's construction (note that part (a) is not immediately
clear: it is conceivable that finite and infinite type geodesics induce the
same orderings).
\begin{theorem}\sl If $S$ is a punctured disk, we have:

(a) An ordering cannot be both of finite {\it and} infinite type.

(b) Given two geodesics $\gamma_\alpha, \gamma_\beta$
of finite type, one can decide whether or not they 
determine the same ordering.

(c) Given two geodesics $\gamma_\alpha, \gamma_\beta$
of infinite type, one can decide whether or not they 
determine the same ordering. For instance, if $\gamma_\alpha$
and $\gamma_\beta$ are embedded, then they determine the same
ordering if and only if 
$\beta=\Delta^{2k}(\alpha)$ for some $k\in \Z$ (ie if $\gamma_\alpha$ and 
$\gamma_\beta$ are related by a slide of the starting point around
$\partial D_n$).
\end{theorem}
In fact, we shall develop machinery which gives a very good and explicit
understanding of finite type orderings:
\begin{theorem} \sl There is only a finite number
of conjugacy classes of orderings of finite type of $\MCG(D_n)=B_n$.
In fact, the number $N_n$ of
conjugacy classes can be calculated by the following recursive formula
$$
N_2=1 \hbox{ \ \ and \ \ } N_{n}=N_{n-1}+\sum_{k=2}^{n-2} 
{{n-2} \choose {k-1}} \ N_k \  N_{n-k}.
$$
\end{theorem}
The following list gives the first few values of $N_n$:

\begin{verbatim}
       n  |    2 |    3 |    4 |    5 |    6 |    7 |    8 |
      N_n |    1 |    1 |    3 |    9 |   39 |  189 | 1197 |
\end{verbatim}

Theorems 3.4 and 3.5 almost certainly generalise to mapping class
groups of other negatively curved surfaces, but in order to keep
our machinery simple, we stick to the special case of punctured disks.


\section{Orderings of mapping class groups using curve diagrams}

In this section we present another method for constructing left-orderings
on $B_n$, using certain diagrams on $D_n$, which we call {\it curve diagrams}.
Both the definition of curve diagrams and the orderings associated to
them are generalisations of similar concepts in \cite{FGRRW}.

{\bf Convention } Whenever we talk about geodesics in $D_n$, we think
of the punctures as being holes in the disk, whose neighbourhoods on
the disk have the geometry of cusps. By contrast, when we talk about
curve diagrams, we think of the punctures as distinguished points 
on, and belonging to, the disk, and we ignore the geometric structure. 
This changing perspective should not cause confusion.

\begin{definition} \sl A (partial) curve diagram $\Gamma$ is a diagram 
on $D_n$ consisting of $j\leqslant n-1$ closed, oriented arcs which are 
labelled $\Gamma_1,\ldots,\Gamma_j$. Moreover, the boundary circle of 
$D_n$ is labelled $\Gamma_0$, and by abuse of notation we shall refer to 
it as an ``arc'' of $\Gamma$. We require: 

(1) every path component of $D_n\basl \Gamma$ has at least one
puncture in its interior, 

(2) $\bigcup_{i=0}^j int(\Gamma_i)$ is embedded and disjoint from the
punctures (where $int$ denotes the interior), 

(3) the starting point of the $i$th arc lies in $\bigcup_{k=0}^{i-1}\Gamma_k$,
ie on one of the previous arcs, 

(4) the end point of the $i$th arc lies in one of the previous arcs, or
on an earlier point of the $i$th arc, or in a puncture. 

In the special case that $j=n-1$, so that in (1) every path component 
contains precisely one puncture, we say $\Gamma$ is a total curve diagram.
\end{definition}
\begin{figure}[htb] %
\cl{
\epsfbox{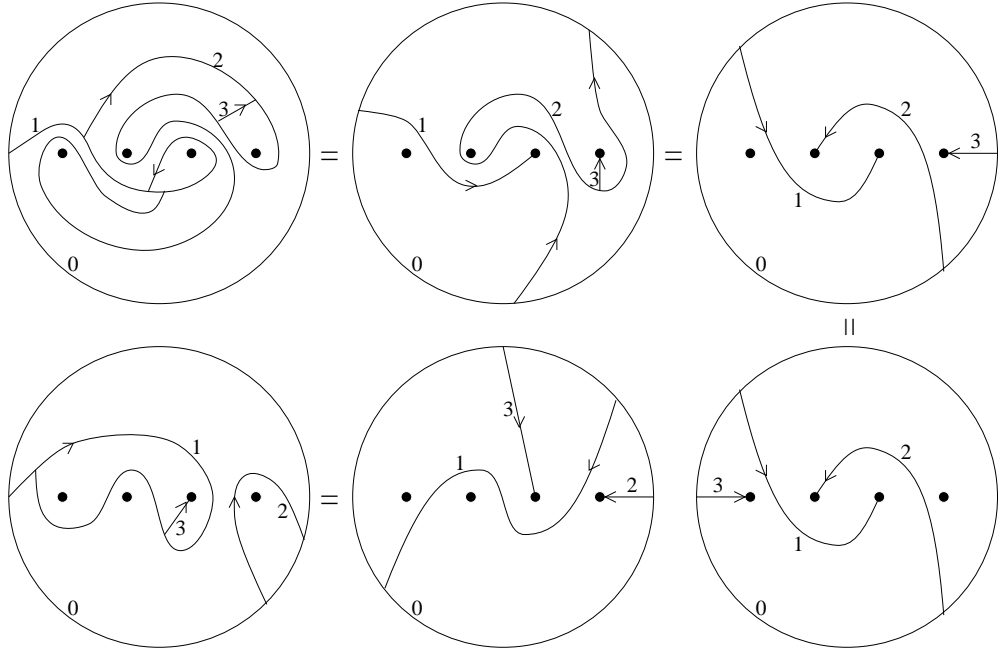} 
}
\caption{Examples of total curve diagrams on $D_4$. The meaning of the 
equality signs will be explained in section 5.}
\label{cdex}
\end{figure}

{\bf Remarks } For simplicity we shall sometimes label arcs
$0,\ldots, j$, instead of $\Gamma_0,\ldots,\Gamma_j$.
Moreover, we shall use the abbreviated notation
$\Gamma_{0\cup\ldots\cup i}:= \bigcup_{k=0}^{i}\Gamma_k$.
Note for (1) that the number of path components of $D_n\basl \Gamma$ equals
$1$ plus the number of arcs of $\Gamma$ not ending in a puncture, so 
it can be anything between $1$ and $n$. Note for (3) that the start
point of the $i$th arc can lie in a puncture, if this puncture was 
the end point of one of the previous arcs. Finally note that if $i<j$ then
$\Gamma_i$ is disjoint from the interior of $\Gamma_j$.

We now explain how to associate a partial left ordering of $\MCG(D_n)=B_n$ 
to a partial curve diagram (with total curve diagrams giving rise to
total orderings).
The essential ingredient in this definition is the well-known procedure
of ``pulling tight'' or ``reducing'' two properly embedded curves in a 
surface. In brief, two simple closed or properly embedded curves in a 
surface can be isotoped into a relative position in which they have minimal 
possible intersection number, and this relative position is unique. Moreover, 
it can be found in a very naive way: whenever one sees a D-disk 
(or ``bigon'') enclosed by a pair of segments of the curves, one ``squashes'' 
it, ie one reduces
the intersection number of the two curves, by isotoping the arcs across
the disk. A systematic exposition of these ideas can for instance be found 
in section 2 of \cite{FGRRW}. 

Our definition of the ordering of $\MCG(S)$ associated to a curve diagram 
will be a (slightly less elegant) variation of the definition in 
\cite{FGRRW}. We briefly remind the reader of this comparison method. 
Let $\Gamma$ be a partial curve diagram in which all arcs are embedded 
(no curve $\Gamma_i$ has end point in its own interior), and let $\phi$ and 
$\psi$ be two homeomorphisms of $D_n$. If $\phi(\Gamma_k)$ and 
$\psi(\Gamma_k)$ are isotopic for $k=1,\ldots, i-1$, then we replace
$\phi$ by an isotopic map, also denoted $\phi$, such that the restrictions
of $\phi$ and $\psi$ to 
$\Gamma_{0\cup\ldots\cup i-1}$ are exactly the same
maps. At this point, $\phi(\Gamma_i)$ and $\psi(\Gamma_i)$ have the
same starting point and lie in the same path component of 
$D_n \basl \phi(\Gamma_{0\cup \ldots \cup i-1})$.
Next we ``pull $\phi(\Gamma_i)$ tight'' 
with respect to $\psi(\Gamma_i)$,
ie\ we isotope $\phi$ so as to minimise the number of intersections
of $\phi(\Gamma_i)$ and $\psi(\Gamma_i)$, as described above.
This can be done by an isotopy which fixes 
$\phi(\Gamma_{0\cup\ldots\cup i-1})$. Restricting finally our attention to
small initial segements of $\phi(\Gamma_i)$ and $\psi(\Gamma_i)$, we
see that the two curves set off from their common starting point into
the interior of a component of $D_n \basl 
\phi(\Gamma_{0\cup \ldots \cup i-1})$ in different directions, one
of them ``going more to the left''; if it is $\phi(\Gamma_i)$ 
say, then we define $\phi>\psi$, otherwise $\psi<\phi$. 
The resulting ordering is left invariant, because the {\it relative}
position of $\chi \circ \phi(\Gamma)$ and $\chi \circ \psi(\Gamma)$ 
is the same as the one of $\phi(\Gamma)$ and $\psi(\Gamma)$ for all 
$\chi \in \MCG(S)$.

\begin{figure}[htb] %
\cl{
\relabelbox
\epsfbox{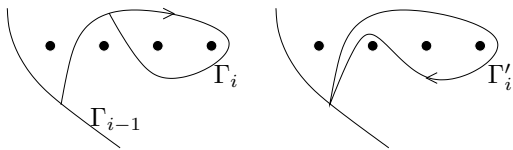} 
\relabel{gammai-1}{\small $\Gamma_{i-1}$}
\relabel{gammai}{\small $\Gamma_i$}
\relabel{gammai'}{\small $\Gamma'_i$}
\endrelabelbox
}
\caption{The embedded arc $\Gamma_i'$ obtained from $\Gamma_i$ by sliding the
end point
}
\label{slback}
\end{figure}

We shall use a variant of this comparison method: first we make 
$\phi(\Gamma_{0\cup\ldots\cup i-1})$ and $\psi(\Gamma_{0\cup\ldots\cup i-1})$
agree for maximal possible $i$, as before.
If the arc $\Gamma_i$ is embedded, we proceed as before, to compare 
$\phi(\Gamma_i)$ and $\psi(\Gamma_i)$.
If the arc $\Gamma_i$ has end point in the interior of $\Gamma_i$ itself,
then we consider the embedded arc $\Gamma_i'$ which, by definition,
is obtained from $\Gamma_i$ by sliding the end point back along $\Gamma_i$
so as to make start and end point coincide, as illustrated in Figure 
\ref{slback}.
We then ignore the original arc $\Gamma_i$, and compare $\phi(\Gamma_i')$
and $\psi(\Gamma_i')$ as before.

\begin{definition} \sl The ordering defined in this way is 
{\it the ordering associated to} the curve diagram $\Gamma$. 
\end{definition}
\begin{lem}\sl The ordering associated to a curve diagram $\Gamma$ is
total if and only if $\Gamma$ is a total curve diagram.
\end{lem}
\begin{proof} If $\Gamma$ is total, ie if all components of 
$D_n \basl \Gamma$ are once-punctured disks, then any homeomorphism of $D_n$
which fixes $\Gamma$ is isotopic to the identity; this follows from the
Alexander trick (see eg \cite{RS}). Conversely, 
if $D_n\basl \Gamma$ has a path component which 
contains at least two holes, then we can push the boundary curve of
this path component slightly into its interior, to make it disjoint
from $\Gamma$. A Dehn twist along such a curve is a nontrivial 
element of $B_n$, and acts trivially on $\Gamma$. 
\end{proof}

{\bf Example } For any $n$, the Dehornoy ordering \cite{D} is defined by
the diagram consisting of $n-1$ horizontal line segments, connecting 
$\partial D_n$ to the first (leftmost) hole, the first to the second hole, 
and so on. The arcs are oriented from left to right, and labelled 
$1,\ldots,n-1$ in this order (see \cite{FGRRW}).

\begin{definition}\sl A (possibly partial) order on a group $G$ is
{\it discrete} if the positive cone $P=\{g\in G \mid g>1\}$ has a minimal 
element. (If the ordering is total then this element is necessarily unique.)
\end{definition}
In a group with a discrete total order every element has a unique predecessor
and successor. We note that an ordering is non-discrete if and only if
for all $a, c \in G$ there exists a $b\in G$ such that $a<b<c$.

\begin{lem} \sl The total ordering associated to a total curve diagram 
$\Gamma$ is discrete. 
\end{lem}
\begin{proof} The curve diagram $\Gamma_{0\cup \ldots \cup n-2}$ (which is 
obtained from $\Gamma$ by removing the arc of maximal index) cuts $D_n$ 
into a number of once-punctured disks and one twice-punctured disk.
We observe that the unique smallest element is the positive half-twist 
interchanging the two punctures inside this disk. \end{proof}

{\bf Remark } It is an easy exercise to prove that the partial orderings 
associated to partial curve diagrams are in general not discrete. However, 
we shall see in the proof of Theorem 3.4(a) that even such orderings have 
a certain discreteness property.


\section{Which pairs of curve diagrams determine the same ordering?}

In this section we define an equivalence relation of
curve diagrams which we call {\it loose isotopy}. 
We give a simple algorithm to decide
whether or not two given curve diagrams are loosely isotopic.
We prove that two curve diagrams determine the same ordering if and
only if they are loosely isotopic. Moreover, the quotient of the set of 
loose isotopy classes of curve diagrams under the natural action of $B_n$
is finite; we deduce that for fixed $n \geqslant 2$ there is 
only a finite number of conjugacy classes of orderings arising from 
curve diagrams.

\begin{definition}\sl Let $\mathcal C$ denote the space of all curve diagrams, 
equipped with the natural topology (the subset topology from the space of 
all mappings of $n-1$ arcs into $D_n$). We define {\it loose isotopy} to 
be the equivalence relation on $\mathcal C$ generated by the following two
types of equivalence:

(1) Continuous deformation: two curve diagrams are equivalent if they
lie in the same path component of $\mathcal C$. 

(2) Pulling loops around punctures tight: if some final segment of the
curve $\Gamma_i$ say cuts out a disk with one puncture from $D_n$,
then this final segment can be pulled tight, so as to make $\Gamma_i$
end in the puncture. 
\end{definition}

\begin{figure}[htb] %
\cl{
\epsfbox{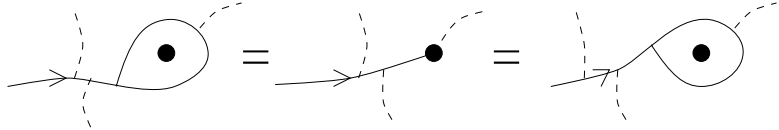} 
}
\caption{Pulling loops around punctures tight}
\label{pupu}
\end{figure}

Equivalence (2) is illustrated in figure \ref{pupu}; here the 
dashed lines indicate any number of arcs of index greater than $i$ which
start on $\Gamma_i$. The reverse of this operation is also allowed.
Equivalence (1) says that one is allowed to
deform the diagram, to slide starting points of arcs along the union
of all previous arcs, including their start and end points, and
even across punctures, if they are the end points of some previous arcs.
Similarly, end points of arcs are allowed to slide across the union
of all ``previous points of the diagram''.

In order to get a feel for the meaning of this definition, the reader may
want to prove that the equality signs in Figure \ref{cdex} represent
loose isotopies.

\begin{theorem}\sl (a) Two curve diagrams determine the same ordering 
of $B_n$ if and only if they are loosely isotopic.

(b) There is an algorithm to decide whether or not two
curve diagrams $\Gamma$ and $\Delta$ are loosely isotopic. \end{theorem}
\begin{proof} For the implication ``$\Leftarrow$'' of (a) we have to prove
that loosely isotopic diagrams define the same ordering. The only 
nonobvious claim here is that the ordering is invariant under the
``pulling tight'' procedure. 

In order to prove this, we consider a curve diagram $\Gamma'$ with $j$
arcs, the $i$th of which is a loop (ie the end point equals the start
point) which encloses exactly one puncture. We consider in addition the
curve diagram $\Gamma$ which is obtained from $\Gamma'$ by squashing
the curve $\Gamma'_i$ to an arc from the starting point of $\Gamma'_i$ to
the enclosed puncture, much as in Figure \ref{pupu}. Let $\phi$ and $\psi$ 
be two nonisotopic homeomorphisms, and more precisely assume that 
$\phi>_\Gamma\psi$. Our aim is to prove that $\phi>_{\Gamma'}\psi$. If 
$\phi(\Gamma_{0\cup\ldots\cup i-1})$ and $\psi(\Gamma_{0\cup\ldots\cup i-1})$
are already nonisotopic then this is obvious since the first $i-1$ arcs
of $\Gamma$ and $\Gamma'$ coincide. On the other hand, if 
$\phi(\Gamma_{0\cup\ldots\cup i})$ and $\psi(\Gamma_{0\cup\ldots\cup i})$
are isotopic (and the difference between $\phi$ and $\psi$ only shows up
on arcs of higher index), then after an isotopy the first $i$ arcs of 
$\phi(\Gamma')$ and $\psi(\Gamma')$ coincide as well, and the result 
follows easily. Finally in the critical case, when the first difference 
occurs on the $i$th arc of $\Gamma$,
we have the two arcs $\phi(\Gamma_i)$ and $\psi(\Gamma_i)$ which are
reduced with respect to each other, with $\phi(\Gamma_i)$ setting off
more to the left. The crucial observation is now that the boundary curves
of sufficiently small regular neighbourhoods of the two curves are isotopic
to $\phi(\Gamma'_i)$ respectively $\psi(\Gamma'_i)$ {\it and} reduced with
respect to each other - see Figure \ref{pulleq}. It is now clear that
$\phi(\Gamma'_i)$ also sets off more to the left than $\psi(\Gamma'_i)$.
This completes the proof of implication ``$\Leftarrow$'' of (a).

\begin{figure}[htb] %
\cl{
\relabelbox
\epsfbox{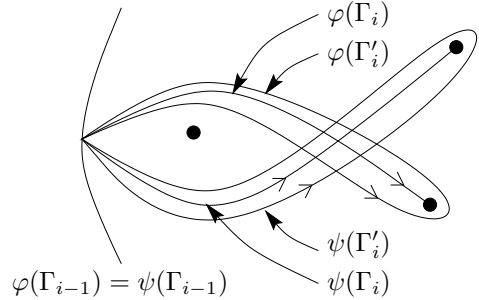} 
\relabel{psii}{\small $\psi(\Gamma_i)$}
\relabel{psii'}{\small $\psi(\Gamma'_i)$}
\relabel{phii}{\small $\phi(\Gamma_i)$}
\relabel{phii'}{\small $\phi(\Gamma'_i)$}
\relabel{phipsii-1}{\small $\phi(\Gamma_{i-1})=\psi(\Gamma_{i-1})$}
\endrelabelbox
}
\caption{Proof that $\phi>_\Gamma\psi \ \Rightarrow \phi>_{\Gamma'}\psi$ - 
the critical case where the first difference between $\phi$ and $\psi$
occurs on the arc which is being pulled tight}
\label{pulleq}
\end{figure}

We shall now explicitely describe the algorithm promised in
(b), and prove the implication ``$\Rightarrow$'' of (a) along the way. 
The proof is by induction on $n$. For the case $n=2$ we note that
any two total curve diagrams (with one arc) are loosely isotopic.
Thus there are only two loose isotopy classes of curve diagrams:
the empty diagram and the one with one arc.
The empty diagram induces the trivial ordering, whereas the diagram with 
one arc induces the ordering $\sigma_1^k > \sigma_1^l \iff k>l$.
So the desired algorithm consists just of counting the number of arcs,
and non loosely isotopic curve diagrams do indeed induce different orderings.

Now suppose that the result is true for disks with fewer than
$n$ punctures, and we want to compare two curve diagrams
$\Gamma_0,\ldots,\Gamma_j$ and $\Delta_0,\ldots,\Delta_{j'}$
in $D_n$, with $j,j' \leq n-1$ and $n\geqslant 3$. 
The arc $\Gamma_1$ ends either on 
$\partial D_n$, or in the interior of $\Gamma_1$ itself, or in a
puncture. In the first two cases $D_n\basl \Gamma_1$ has precisely
two path components. At most one of them can contain only one puncture;
if one of them does, we pull $\Gamma_1$ tight around it. If both
components of $D_n \basl \Gamma_1$ contain more than one puncture and if 
$\Gamma_1$ ends on itself, then we slide the end point of $\Gamma_1$ back 
along $\Gamma_1$, across its starting point, and into $\Gamma_0=\partial 
D_n$. There are now two possibilities left: either $\Gamma_1$ is an embedded 
arc connecting the boundary to a puncture ($\Gamma_1$ is {\it nonseperating}), 
or it is an embedded arc connecting two boundary points, cutting $D_n$ 
into two pieces, each of which has at least two punctures in its interior 
($\Gamma_1$ is {\it seperating}). We repeat this procedure for $\Delta_1$. 
There are now four cases:

(1) It may be that $\Gamma_1$ is seperating, while $\Delta_1$ is not 
(or vice versa).

(2) It is possible that $\Gamma_1$ and $\Delta_1$ are both nonseperating 
but are not isotopic with starting points sliding in $\partial D_n$ 
(a criterion which is easy to check algorithmically). 

(3) It is possible that $\Gamma_1$ and $\Delta_1$ are both seperating but 
are not isotopic as oriented arcs, with starting and end points sliding 
in $\partial D_n$ (a criterion which is equally easy to check 
algorithmically). 

{\bf Claim } {\sl In these first three cases the orderings defined by 
$\Gamma$ and $\Delta$ do not coincide, and $\Gamma$ and $\Delta$ are
not loosely isotopic.}

We only need to prove the first part of the claim, the second one follows
by the implication ``$\Leftarrow$'' of Theorem 5.2(a).
We first treat the following pathological
situation: if, in case (3) above, $\Gamma_1$ and $\Delta_1$ are isotopic 
to each other, but with opposite orientations, then a homeomorphism of 
the type indicated in Figure \ref{ga1de1} 
\begin{figure}[htb] %
\cl{
\relabelbox
\epsfbox{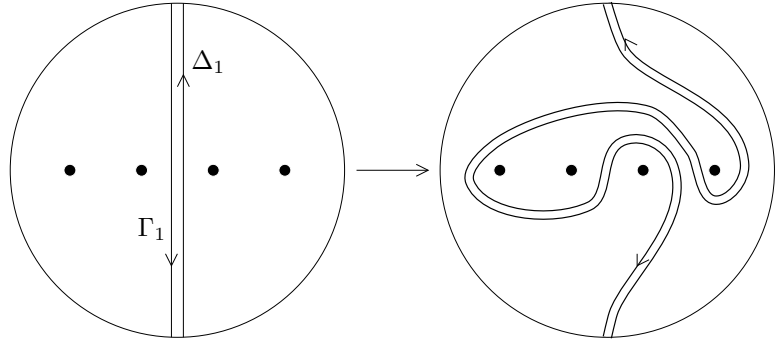} 
\relabel{de1}{\small $\Gamma_1$}
\relabel{ga1}{\small $\Delta_1$}
\endrelabelbox
}
\caption{A homeomorphism which distinguishes the $\Gamma$- and 
$\Delta$-orderings}
\label{ga1de1}
\end{figure}
is positive in the ordering defined by $\Gamma$,
but negative in the $\Delta$-ordering. In all other situations allowed by 
(1), (2), and (3), there exists a simple closed curve $\tau$ in $D_n$ which is 
disjoint from $\Gamma_1$, but intersects every arc isotopic to $\Delta_1$.
(Consider, for instance, a regular neighbourhood of 
$\partial D_n \cup \Gamma_1$ in $D_n$. If $\Gamma_1$ is nonseperating then 
its boundary curve has this property; if $\Gamma_1$ is seperating then at 
least one of the two boundary curves has.) 
We denote by $T \co D_n \to D_n$ the positive Dehn twist along $\tau$.
The map $T$ leaves $\Gamma_1$ invariant, while the arc
$T(\Delta_1)$ is ``more to the left'' than $\Delta_1$ by Proposition 2.4.
Similarly, there exists a curve $\tau'$ which is disjoint from 
$\Delta_1$, but not from any arc in the isotopy class of $\Gamma_1$. 
Then ${T'}\inv$ sends $T(\Delta_1)$ more to the right, but not very far:
${T'}\inv \circ T(\Delta_1)$ is still to the left of 
the arc $\Delta_1$, which is fixed by ${T'}\inv$; and ${T'}\inv$ sends 
$\Gamma_1$ to the right, as well. Thus, in summary, the composition
${T'}\inv \circ T$ sends $\Delta_1$ more
to the left but $\Gamma_1$ more to the right, so that
${T'}\inv\circ T \in B_n$ is negative in the ordering determined
by $\Gamma$, but positive in the $\Delta$-ordering. This proves the claim.
(One may find simpler proofs, but this one will be useful in section 7.)

(4) The remaining possibility is that $\Gamma_1$ and $\Delta_1$ 
{\it can} be made to coincide by isotopies which need not be fixed on 
$\partial D_n$. Such isotopies can be extended to loose isotopies of 
$\Gamma$ or $\Delta$.

To summarize, we can algorithmically decide whether or not there
is a loose isotopy which makes $\Gamma_1$ and $\Delta_1$ coincide.
If the answer is No (cases (1) - (3)), then $\Gamma$ and $\Delta$ are 
not loosely isotopic, and the orderings defined by $\Gamma$ and 
$\Delta$ do not coincide. In this case, the implication ``$\Rightarrow$'' 
of 5.2(a) is true. If the answer is Yes (case (4)), then 
$D_n\basl \Gamma_1 = D_n\basl \Delta_1$ has either one or two path 
components, each of which is a disk with at most $n-1$ punctures. 
Moreover, the arcs $\Gamma_2,\ldots,\Gamma_j$ form curve diagrams 
in these disks (with some indices missing in each curve diagram, if the 
arcs are distributed among two disks),
and similarly for $\Delta_2,\ldots,\Delta_j$.
Finally, the following conditions are equivalent:\\
(i) $\Gamma$ and $\Delta$ are loosely isotopic, \\
(ii) in each path component of $D_n\basl \Gamma_1 = D_n\basl \Delta_1$
there is a loose isotopy between the diagrams made up of the remaining
arcs of $\Gamma$ respectively $\Delta$,\\
(iii) the orderings of $Fix(\Gamma_1) \subseteq B_n$ induced by $\Gamma$
and $\Delta$ coincide, where $Fix(\Gamma_1)$ denotes the subgroup whose
elements have support disjoint from $\Gamma_1$,\\  
(iv) the orderings of $B_n$ defined by $\Gamma$ and $\Delta$ coincide.\\
The equivalences between (i) and (ii), and between (iii) and (iv) are
clear, whereas the equivalence of (ii) and (iii) follows from the induction 
hypothesis. Also by induction hypothesis, we can decide algorithmically 
whether or not (ii) holds. This proves the theorem in case (4).
\end{proof}

We recall that for any ordering ``$<$'' of $B_n$, and every element 
$\rho \in B_n=\MCG(D_n)$ one can construct an ordering ``$<_\rho$'', 
by defining 
$\phi<_\rho \psi$ $:\hspace{-1mm}\iff$ $\phi \rho < \psi \rho$, 
and we call $<_\rho$ ``the ordering $<$ conjugated by $\rho$''. We observe 
that if $<$ is induced by a curve diagram $\Gamma$, then $<_\rho$
is induced by the curve diagram $\rho(\Gamma)$. Thus two curve diagrams
$\Gamma$ and $\Delta$ induce conjugate orderings if and only if 
$\Gamma$ and $\Delta$ are in the same orbit under the natural action of
$B_n$ on the set of loose isotopy classes of curve diagrams.

\begin{prop}\sl Let $M_n$ denote the number of conjugacy classes of total 
orderings of $B_n$ arising from curve diagrams. Then $M_n$ can be 
calculated by the following recursive formula
$$
M_2=1 \hbox{ \ \ and \ \ } M_{n}=M_{n-1}+\sum_{k=2}^{n-2} 
{{n-2} \choose {k-1}} \ M_k \  M_{n-k}.
$$
\end{prop}
{\bf Remark } In order to avoid confusion, we recall our orientation 
convention: we are insisting that ``more to the left'' means ``larger''.
It is for this reason that there is only one ordering of $B_2=\Z$, not
two, as one might expect.
\begin{proof} We shall count the orbits of
the set of loose isotopy classes of total curve diagrams under the action
of $B_n$. The case $n=2$ is clear, since there is only one loose
isotopy class of curve diagrams. Now suppose inductively that the
formula is true for up to $n-1$ strings. 

For every total curve diagram in $D_n$ there are two possibilities: 
(a) the first arc of the curve diagram ends in a puncture or can be 
pulled tight so as to end in a puncture; 
(b) the first arc cuts $D_n$ into two disks, each of which contains at 
least two punctures.

For case (a) we notice that the first arc 
can be turned into the horizontal arc from $-1$ to the leftmost
puncture, by an action of some appropriate element of $B_n$. There are now
precisely $M_{n-1}$ orbits of loose isotopy classes of curve diagrams of 
the remaining $n-2$ arcs in the $n-1$-punctured disk 
$D_n \basl$(the first arc). So case (a) gives a contribution of $M_{n-1}$
orbits.

The argument for case (b) is similar: the action of an appropriate element
of $B_n$ will turn the first arc of any curve diagram of type (b) into
the vertical arc, oriented from bottom to top, having $k$ punctures on its
left and $n-k$ on its right, for some $k \in \{2,\ldots , n-2\}$. 
In this case, there should be $k-1$ arcs on the left and $n-k-1$ arcs on the 
right of the first arc, so there are ${{n-2} \choose {k-1}}$ ways to 
distribute the remaining $n-2$ arcs over the two sides. Finally, there
are $M_k$ respectively $M_{n-k}$ orbits of loose isotopy classes of total 
curve diagrams on the disk on the left respectively on the right.
\end{proof}


\section{Replacing geodesics by curve diagrams}

In this section we prove the main theorems on orderings of finite type. 
The strategy is to associate  to every geodesic 
of finite type a curve diagram such that the (possibly partial) orderings 
arising from  the geodesic and the curve diagram agree. Thus we obtain, 
via curve diagram orderings, a good understanding of finite type orderings.

\begin{proof}[Proof of Theorem 3.3(a)]
If $D_n\basl \gamma_\alpha$ has a path component which 
contains at least two holes, then we can push the boundary curve of
this path component slightly into its interior, to make it disjoint
from $\gamma_\alpha$. A Dehn twist along such a curve will be a nontrivial 
element of $B_n$, and act trivially on $\gamma_\alpha$. 
\end{proof}

We now define {\it the curve diagram $C(\gamma_\alpha)$ associated
to a geodesic} $\gamma_\alpha$ of finite type.
It is a subset of $\gamma_\alpha$, more precisely a union of segments of 
$\gamma$ which start and end at self-intersection points. The diagram will 
be disjoint from the punctures, except that the last arc may fall into a 
puncture. For simplicity we shall write $\Gamma$ for $C(\gamma_\alpha)$
and, as before, $\Gamma_{0 \cup \ldots \cup i-1}$ for 
$\bigcup_{k=0}^{i-1} \Gamma_k$,

\begin{figure}[htb] %
\cl{
\epsfbox{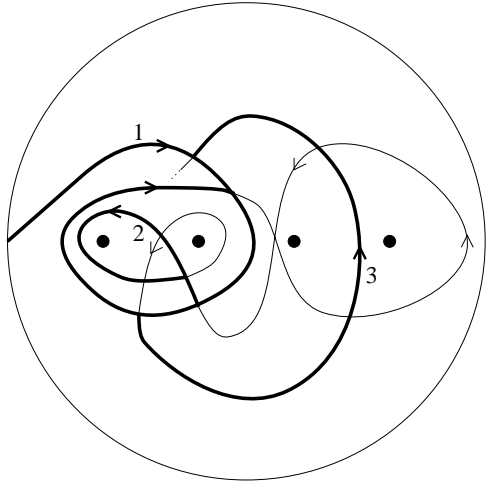} 
}
\caption{A geodesic and (in bold line) its associated curve diagram}
\label{geotocd}
\end{figure}

The definition is inductive. We define $\Gamma_0 = \partial D_n$. Now 
suppose that we have already found $\Gamma_0,\ldots ,\Gamma_{i-1}$.
So every path component of
$D_n \basl \Gamma_{0 \cup \ldots \cup i-1}$ is a disk containing at least
one puncture. We put down a pencil at the end point of $\Gamma_{i-1}$, 
start tracing out $\gamma_\alpha$, drawing an arc $\Gamma_i^p$ (with 
``$p$'' standing for ``potential'', because $\Gamma_i^p$ is potentially the
new arc $\Gamma_i$). We continue drawing either up to the next intersection
with $\Gamma_{0 \cup \ldots \cup i-1}$, or up to the first self intersection
of $\Gamma_i^p$, or until $\gamma_\alpha$ falls into a puncture, whichever 
comes first. We now decide whether or not
$\Gamma_i^p$ has cut one of the components of $D_n \basl 
\Gamma_{0\cup\ldots\cup i-1}$ in a nontrivial way, 
ie whether it has either fallen into a puncture or cut one
of the components of $D_n \basl \Gamma_{0 \cup \ldots \cup i-1}$ into
two, both of which contain at least one puncture. If yes, we let
$\Gamma_i:= \Gamma_i^p$, and have finished the induction step. If not,
we rub out $\Gamma_i^p$, and start a new $\Gamma_i^p$ at the
next intersection point of $\gamma_\alpha$ with 
$D_n \basl \Gamma_{0 \cup \ldots \cup i-1}$. (This intersection point is just 
the end point of the previous $\Gamma_i^p$, unless the previous $\Gamma_i^p$ 
had end point in its own interior. 
Note that in this latter case not only $\Gamma_i^p$, but the entire segment 
of the geodesic $\gamma_\alpha$ between successive intersection points with
$\Gamma_{0 \cup \ldots \cup i-1}$ cuts the disk in a trivial way).

\begin{figure}[htb] %
\cl{
\relabelbox
\epsfbox{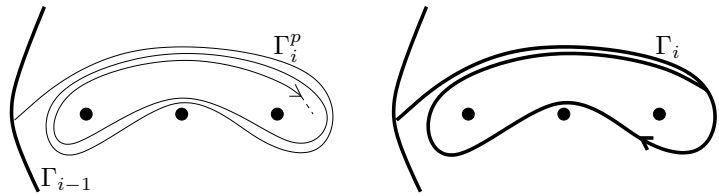} 
\relabel{Gammaip}{\small $\Gamma_i^p$}
\relabel{Gammai}{\small $\Gamma_i$}
\relabel{Gammai-1}{\small $\Gamma_{i-1}$}
\endrelabelbox
}
\caption{The curve diagram associated to a geodesic which spirals towards
a closed geodesic}
\label{spiral}
\end{figure}

There is one special rule: if in the construction process we obtain an
arc $\Gamma_i^p$ which spirals {\it ad infinitum} towards a simple closed
geodesic, then we define $\Gamma_i$ to be the arc with end point in
its own interior containing $\Gamma_i^p$ in a regular neighbourhood,
as shown in figure \ref{spiral} (this arc is unique up to loose isotopy).

Since at most $n-1$ arcs can be constructed in this way, the process
terminates after finitely many steps.
We observe that the curve diagram $C(\gamma_\alpha)$ is total if and only if
the geodesic $\gamma_\alpha$ fills $D_n$. More generally, two punctures are
in the same path component of $D_n \basl \gamma_\alpha$ if and only if
they are in the same path component of $D_n \basl C(\gamma_\alpha)$.
We also note that for every geodesic $\gamma_\alpha$ and $\phi \in B_n$ we 
have $C(\phi(\gamma_\alpha))=\phi(C(\gamma_\alpha))$.

\begin{theorem}\sl For any $\alpha \in (0,\pi)$ and $\phi \in B_n$ we have:

(a) if the curve diagrams $\phi(C(\gamma_\alpha))$ and 
$C(\gamma_\alpha)$ are isotopic then $\phi(\alpha)=\alpha$.

(b) if $\phi(C(\gamma_\alpha)) > C(\gamma_\alpha)$ 
(in the curve diagram sense) then we have $\phi(\alpha)>\alpha$ in $\R$.
\end{theorem}
\begin{cor}\sl For every geodesic $\gamma_\alpha$ of finite type 
($\alpha \in (0,\pi)$), the ordering of $B_n$ associated to $\alpha$ by 
Remark 1.2(1) coincides with the ordering associated to the curve diagram 
$C(\gamma_\alpha)$ by Definition 4.2. \end{cor}

\begin{proof}[Proof of the theorem] 
We shall need a generalisation of the concept of relative ``reduction'' of
two simple curves in $D_n$ to the case where one of the two curves is
authorised to have self-intersections, but no D-disks with itself.
For instance, we shall be interested in the case where one of the two
curves is a simple geodesic, and the other is a homeomorphic image of a 
non-simple geodesic.

Suppose that $C$ is a disjoint collection of simple closed geodesics
and properly embedded geodesic arcs connecting distinct punctures 
in $D_n$. Then we say that $\phi(\gamma_\alpha)$ is {\it reducible} 
with respect to $C$ if there are D-disks enclosed by $\phi(\gamma_\alpha)$
and $C$, ie if there are finite segments of $\phi(\gamma_\alpha)$
and of $C$ with the same start and end points which are homotopic
with fixed end points. If $\phi(\gamma_\alpha)$ is not reducible
then we say it is {\it reduced} with respect to $C$. Equivalently,
any component of the preimage of $\phi(\gamma_\alpha)$ in the universal 
cover $D_n^\sim$ intersects any component of the preimage of $C$  
at most once.

\begin{lem}\sl One can pull $\phi(\gamma_\alpha)$ tight with respect
to $C$, ie there exists an isotopy of $\phi$ which makes $\phi(\gamma_\alpha)$
and $C$ reduced with respect to each other.
\end{lem}
\begin{proof} The proof is an easy exercise - it is in fact similar to
the proof of the ``triple reduction lemma'' 2.1 of \cite{FGRRW}.
\end{proof}

We need some more notation. We still write $\Gamma$ for $C(\gamma_\alpha)$,
denote by $j$ the number of arcs of $\Gamma$,
and consider the partial curve diagrams $\Gamma_{0 \cup \ldots \cup i-1}$ 
for $i\in \{1,\ldots,j\}$; all their arcs are geodesics. Every path
component of $D_n\basl \Gamma_{0\cup \ldots \cup i-1}$ contains at least 
one puncture in its interior. The boundary curve of each component with at
least two punctures is isotopic to a unique simple closed geodesic, which
bounds a disk (with these punctures in its interior) in $D_n$.
Removing all these disks from $D_n$ yields a planar surface with a 
number of geodesic boundary components (one of them being $\partial D_n$,
the others corresponding to the at least twice punctured components of 
$D_n\basl \Gamma_{0\cup \ldots \cup i-1}$) and a number of punctures
(corresponding to once-punctured components of 
$D_n\basl \Gamma_{0\cup \ldots \cup i-1}$). We denote this surface by
$N\Gamma_{0\cup \ldots \cup i-1}$; it is a regular neighbourhood of 
$\partial D_n \cup \Gamma_{0\cup \ldots \cup i-1}$ in $D_n$, and 
contains the complete initial segment of the geodesic $\gamma_\alpha$
up to the starting point of the arc $\Gamma_i \subset \gamma_\alpha$.

We are now ready to prove the theorem. For part (a) suppose that we are
given $\alpha\in (0,\pi)$, and $\phi\in B_n$, and that the curve diagrams 
$\Gamma$ and $\phi(\Gamma)$ are isotopic. Then we can modify the map 
$\phi$ by an isotopy which fixes $\partial D_n$ such that the restriction 
$\phi|N\Gamma$ becomes the identity map. But by construction of 
$\Gamma= C(\gamma_\alpha)$, the geodesic $\gamma_\alpha$ is entirely
contained in $N\Gamma$, and is thus mapped identically. This proves
part (a) of the theorem.

For part (b) suppose that we are given $\alpha\in (0,\pi)$ and $\phi\in B_n$, 
and that for some $i\in \{1,\ldots,j\}$ the curve diagrams 
$\Gamma_{0\cup \ldots \cup i-1}$ and $\phi(\Gamma_{0\cup \ldots \cup i-1})$
are isotopic,  whereas $\phi(\Gamma_i)$ is ``more to the left'' than 
$\Gamma_i$.
Our aim is to prove that $\phi(\alpha)>\alpha$, ie that the end points
of the liftings of $\phi(\gamma_\alpha)$ and $\gamma_\alpha$ on
$\partial D_n^\eqsim \basl \Pi \cong (0,\pi)$ are different, with the one
of $\phi(\gamma_\alpha)$ being ``higher'' in Figure \ref{hyp}. 

Firstly, the map $\phi$ sends $\Gamma_{0\cup \ldots \cup i-1}$ to
a curve diagram which is isotopic to $\Gamma_{0\cup \ldots \cup i-1}$;
therefore we can assume, after an isotopy of $\phi$ which fixes 
$\partial D_n$, that the restriction $\phi|N\Gamma_{0\cup \ldots \cup i-1}$ 
is the identity map. Note that $\gamma_\alpha$, being a geodesic, is 
already reduced with respect to the collection of 
geodesics $\partial N\Gamma_{0\cup \ldots \cup i-1}$, 
and therefore $\phi(\gamma_\alpha)$ is also reduced with respect to 
$\partial N\Gamma_{0\cup \ldots \cup i-1}=
\phi(\partial N\Gamma_{0\cup \ldots \cup i-1})$.

Next, we note that the arc $\Gamma_i$ will cut precisely one of
the components of $D_n \basl N\Gamma_{0\cup \ldots \cup i-1}$ in two,
and leave the other components untouched. This critical component is
an at least twice punctured disk, and we shall denote it by $D_c$.
The preimage of $D_c$ in the universal cover $D_n^\sim$ has many
path components, but we shall be interested in one particular component
$D_c^\sim$, namely the one which is cut in two by the
segment corresponding to $\Gamma_i\subset \gamma_\alpha$ in the geodesic
$\widetilde{\gamma}_\alpha$ in $D_n^\sim$. 

We now distinguish three cases: firstly, the arc $\Gamma_i$ 
falls into a puncture inside $D_c$;
secondly, the arc $\Gamma_i$ has its end point in 
$N\Gamma_{0\cup\ldots\cup i-1}$ (either on $\Gamma_{0\cup\ldots\cup i-1}$
or in the initial segment $\Gamma_i\cap N\Gamma_{0\cup\ldots\cup i-1}$ 
of $\Gamma_i$); thirdly, the end point of the arc $\Gamma_i$ lies in the 
interior of $D_c$ (and then necessarily in the interior of $\Gamma_i$).

The first case is the easiest: by an isotopy of $\phi$ which is
fixed outside $D_c$ we can pull $\phi(\Gamma_i)\cap D_c$ tight with
respect to $\Gamma_i\cap D_c$. The effect of this isotopy is to make
the images of the liftings $\widetilde{\phi}(\widetilde{\gamma}_\alpha)
\cap \widetilde{D}_c$ and $\widetilde{\gamma}_\alpha \cap \widetilde{D}_c$ 
disjoint, except for the common starting point. Moreover,
both liftings run inside $\widetilde{D}_c$ all the way to the circle
at infinity. By the hypothesis that $\phi(\Gamma)>\Gamma$ we have that an 
initial segment of 
$\widetilde{\phi}(\widetilde{\gamma}_\alpha)$ lies to the left of the 
corresponding segment $\widetilde{\gamma}_\alpha$, we conclude that
its end point on the circle at infinity also lies more to the left.
This proves the theorem in the first case. 

\begin{lem}\sl 
If $\gamma$ is a (finite or infinite) geodesic starting on the boundary 
of the punctured disk $D_c$, and if $\phi$ is an automorphism of $D_c$
which acts nontrivially on $\gamma$,
then two liftings of $\gamma$ and $\phi(\gamma)$
to the universal cover $D_c^\sim$ of $D_c$ with the same starting 
point in $\partial D_c^\sim$ have end points either on different components 
of $\partial D_c^\sim$ (if $\gamma$ is finite) or on different points
at infinity (if $\gamma$ is infinite). \qed
\end{lem}

In the second case, we can pull the arc $\phi(\Gamma_i)\cap D_c$ tight
with respect to $\Gamma_i\cap D_c$ by an isotopy of $\phi$ as in the first
case, thus making their liftings disjoint. We now have by hypothesis
that the point of intersection of $\widetilde{\phi}(\Gamma_i^\sim)$ with
$\partial D_c^\eqsim$ where $\widetilde{\phi}(\Gamma_i^\sim)$ exits
$D_c^\eqsim$ lies to the left of the one of $\Gamma_i^\sim$.
By the previous lemma, the two points will even lie on different 
boundary components of $D_c^\sim$, and therefore there is a point of
$\partial D_c^\eqsim$ between these two boundary components which lies
on the circle at infinity. For the liftings of our geodesic and its image
this means the following: $\widetilde{\gamma}_\alpha$ and 
$\widetilde{\phi}(\widetilde{\gamma}_\alpha)$
enter $\partial D_c^\sim$ at the same point, but exit into different 
components of $D_n^\sim \basl D_c^\sim$, with 
$\widetilde{\phi}(\widetilde{\gamma}_\alpha)$
choosing the one that lies more to the left. Since 
$\widetilde{\gamma}_\alpha$ and $\widetilde{\phi}(\widetilde{\gamma}_\alpha)$
intersect $\partial D_c^\sim$
only once, they stay inside their chosen component
of $D_n^\sim \basl D_c^\sim$. Hence we have for their end points that 
$\phi(\alpha)>\alpha$, and the theorem is proved in the second case.

In the third case we consider the arc $\Sigma :=\Gamma_i'$ as in figure 
\ref{slback},
and for simplicity we choose $\Sigma$ to be a geodesic arc. We denote
by $D_{cc}\subset D_c$ the subdisk cut off by $\Sigma$ (so that
$\Sigma=\partial D_{cc}$). Since $\Sigma$ is geodesic, we have that 
$\gamma_\alpha \cap D_c$ is reduced with respect to $\Sigma$.
After an isotopy of $\phi$ inside $D_c$ we can assume 
by lemma 6.3 that the first component of $\phi(\gamma_\alpha)\cap D_c$ 
(the one that contains $\phi(\Gamma_i)\cap D_c$) is also reduced with respect 
to $\Sigma$. By the hypothesis that $\phi(\Gamma_i)$ sets off more to the 
left than $\Gamma_i$,
we are now in one of the situations indicated in figure \ref{DcDcc}. 

\begin{figure}[htb] %
\cl{
\relabelbox
\epsfbox{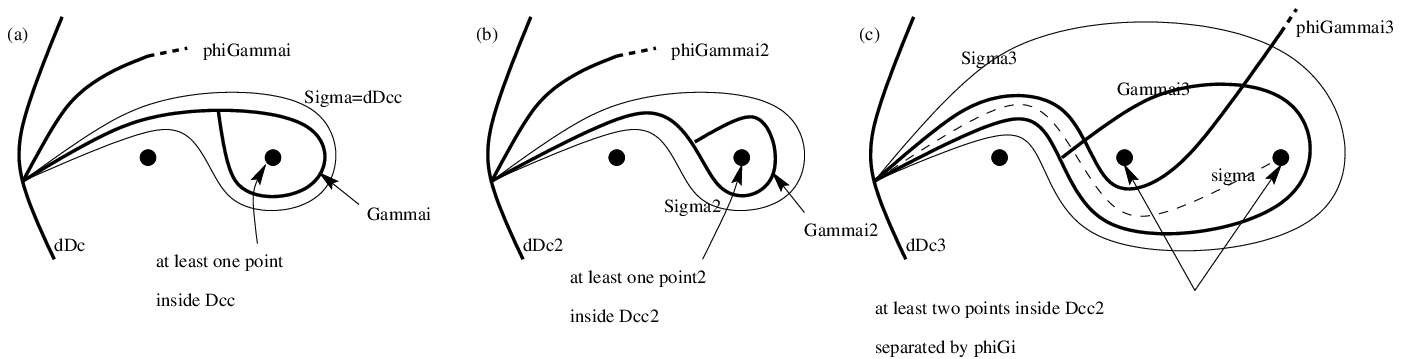} 
\relabel{(a)}{\small (a)}
\relabel{dDc}{\small $\partial D_c$}
\relabel{at least one point}{\small at least one puncture} 
\relabel{inside Dcc}{\small inside $D_{cc}$}
\relabel{Gammai}{\small $\Gamma_i$}
\relabel{Sigma=dDcc}{\small $\Sigma=\partial D_{cc}$}
\relabel{phiGammai}{\small $\phi(\Gamma_i)$}
\relabel{(b)}{\small (b)}
\relabel{phiGammai2}{\small $\phi(\Gamma_i)$}
\relabel{Sigma2}{\small $\Sigma$}
\relabel{at least one point2}{\small at least one puncture} 
\relabel{inside Dcc2}{\small inside $D_{cc}$}
\relabel{dDc2}{\small $\partial D_c$}
\relabel{Gammai2}{\small $\Gamma_i$}
\relabel{(c)}{\small (c)}
\relabel{phiGammai3}{\small $\phi(\Gamma_i)$}
\relabel{sigma}{\small $\sigma$}
\relabel{Sigma3}{\small $\Sigma$}
\relabel{Gammai3}{\small $\Gamma_i$}
\relabel{at least two points inside Dcc2}{\small at least two punctures inside}
\relabel{separated by phiGi}{\small  $D_{cc}$, separated by \thinspace 
$\phi(\Gamma_i)$}
\relabel{dDc3}{\small $\partial D_c$}
\endrelabelbox
}
\caption{The critical disk $D_c$ containing $\Gamma_i$ and $\phi(\Gamma_i)$}
\label{DcDcc}
\end{figure}

A first possiblity
is that an initial segment of $\phi(\Gamma_i)\cap D_c$ lies to the left of 
the tip of $D_{cc}$ (figure \ref{DcDcc}(a) and (b)); in the universal cover 
$D_c^\sim$ we now have three arcs, namely 
$\widetilde{\phi}(\widetilde{\gamma}_\alpha)\cap D_c^\sim$, a lifting of 
$\Sigma$, and $\widetilde{\gamma}_\alpha\cap D_c^\sim$ (and, in fact,
a fourth arc, another lifting of $\Sigma$) starting at the same point 
of $\partial D_c^\sim$, and setting off into different directions, namely
in the given order from left to right. Moreover, the liftings of $\Sigma$
are disjoint from the interiors of the other two arcs, by reducedness. Thus 
the end point of $\widetilde{\phi}(\widetilde{\gamma}_\alpha)\cap D_c^\sim$ 
on $\partial D_c^\eqsim$ lies more to the left than the one of 
$\widetilde{\gamma}_\alpha\cap D_c^\sim$. Even stronger, by Lemma 6.4
they lie either on different points at infinity (in which case we are done)
or they leave $D_c^\sim$ through different components of $\partial D_c^\sim$
(in which case we argue as above that their remainders are trapped in 
different components of $D_n^\sim \basl D_c^\sim$, so that 
$\widetilde{\phi}(\widetilde{\gamma}_\alpha)$ stays to the left of 
$\widetilde{\gamma}_\alpha$). 

The second possibility is that some initial segment of 
$\phi(\Gamma_i)\cap D_c$ lies in $D_{cc}$ (figure \ref{DcDcc}(c)); 
then $D_{cc}$, cut along this 
initial segment, has precisely two path components, each of which contains 
at least one puncture. Since $\phi(\Gamma_i)$ is oriented, we can refer to 
them as the ``left'' and the ``right'' half of $D_{cc}$.
We now consider a geodesic arc $\sigma$ which is embedded in the right half of
$D_{cc}$, starts at the tip of $D_{cc}$ (ie at the same point as 
$\Gamma_i\cap D_c$ and $\phi(\Gamma_i)\cap D_c$), and falls into one of the 
punctures in the right half of $D_{cc}$. By construction, $\gamma_\alpha
\cap D_{cc}$ is reduced with respect to $\sigma$, since both are geodesics,
and the first component of $\phi(\gamma_\alpha)\cap D_{cc}$ is even disjoint
from $\sigma$.
In the universal cover we now have that the lifting $\widetilde{\sigma}$
of $\sigma$ ends on the circle at infinity, thus separating 
$\widetilde{D}_{cc}$ into two components, the left one containing the
lift of $\phi(\gamma_\alpha)\cap D_{cc}$, and the right one the lift
of $\gamma_\alpha\cap D_{cc}$. (Note that this holds in particular
when $\Gamma_i$ represents an infinite geodesic spiralling towards the
simple closed geodesic homotopic to $\Sigma$.)
Thus lifts of these two curves, 
not being allowed to intersect any component of 
$\partial D_{cc}^\sim$ and $\partial D_c^\sim$ more than once, 
go on to hit different points of $\partial D_n^\sim$, with
$\widetilde{\phi}(\widetilde{\gamma}_\alpha)$ staying more to the left
than $\widetilde{\gamma}_\alpha$. This completes the proof of the theorem
in the third case.
\end{proof}

\begin{proof}[Proof of Theorem 3.3(b)] 
If $\gamma_\alpha$ fills $D_n$,
then $C(\gamma_\alpha)$ is a total curve diagram, and thus induces a
{\it total} ordering of $B_n$. By corollary 6.2, the ordering of $B_n$
associated to the point $\alpha \in (0,\pi)$ agrees with this ordering.
\end{proof}
\begin{proof}[Proof of Theorem 3.4(b)] 
For any two geodesics $\gamma_\alpha$ 
and $\gamma_\beta$ of finite type 
one can work out their
associated curve diagrams $C(\gamma_\alpha)$ and $C(\gamma_\beta)$.
By corollary 6.2 it is sufficient to decide whether or not the
orderings associated to the two curve diagrams coincide, which can
be done by Theorem 5.2. 
\end{proof}
\begin{proof}[Proof of Theorem 3.5] It only remains to be proved that 
$N_n=M_n$ (where $M_n$ is given in Theorem 5.3), ie that 
every curve diagram is realized up to loose isotopy as $C(\gamma_\alpha)$ 
for some geodesic $\gamma_\alpha$, $\alpha \in (0,\pi)$. This is left
as an exercise to the reader. 
\end{proof}


\section{Orderings associated to geodesics of infinite type}

In this section we prove the results concerning orderings of infinite
type, and explain the essential differences between finite and
infinite type orderings.

We start by describing in more detail than in section 3 the structure
of geodesics of infinite type. We define the {\it curve diagram 
$C(\gamma_\alpha)$ associated to a geodesic of infinite type}
by precisely the same inductive construction procedure as in the
finite type case. Except for a finite initial segment, the last arc 
$\Gamma_j$ will lie in some path component $D_c$ of 
$D_n \basl N\Gamma_{0\cup \ldots \cup j-1}$, the only difference with
the finite type case is that $\Gamma_j$ goes on for ever, without
falling into a puncture and without spiralling. 
The closure of $\Gamma_j$ inside this critical component $D_c$ is
a geodesic lamination; the lamination has no closed leaves, for such
a leaf would have to be the limit of an infinite spiral of $\Gamma_j$
(see \cite{CB}).
All self-intersections of the geodesic $\gamma_\alpha$ occur inside
the finite initial segment up to the entry into the punctured disk 
$D_c$; in particular, there is only a finite number of self-intersections.
\begin{proof}[Proof of Theorem 3.3(c)] We are studying the set
$${\mathcal I}:= \{ \alpha \in (0, \pi) \| \gamma_\alpha 
\hbox{ \rm  is of infinite  type } \}.$$
 The proof uses standard results
from the theory of geodesic laminations and the Nielsen-Thurston 
classification of surface automorphisms \cite{CB}. 

That $\mathcal{I}$ has uncountably many elements
follows from the fact that there are
uncountably many geodesic laminations of $D_n$, only countably many
of which fall into infinite spirals. A more practical way of seeing
this is to choose arbitrarily a fundamental domain of $D_n$ by
fixing $n$ geodesic arcs, as eg shown in figure \ref{hyp}. 
Thus the fundamental domain is a $2n+1$-gon with one boundary edge
corresponding to $\partial D_n$ and $n$ pairs of boundary edges which
are identified in $D_n$. A segment of the geodesic between any two
sucessive intersections with the boundary of the fundamental domain 
consists of an embedded arc connecting different edges of the $2n+1$-gon.
Hence constructing a geodesic of infinite type amounts to 
choosing an infinite ``cutting sequence'' of the geodesic with 
the boundary arcs of the fundamental domain. Often the
choice will be forced upon us by the requirement that the geodesic
be embedded, but there will be an infinite number of times when we
have a genuine choice. Thus the set of all possible sequences of 
choices is uncountable.

The cutting sequence approach also makes it clear why 
${\mathcal I} \subset (0,\pi)$ is discrete and its elements are nonisolated. 
Given an $\epsilon>0$, there exists an $N_\epsilon\in \N$ such that all 
geodesics $\gamma_\beta$ whose cutting sequences agree with the one of 
$\gamma_\alpha$ for at least $N_\epsilon$ terms satisfy 
$|\alpha - \beta| < \epsilon$. 
Now elements of $\mathcal{I}$ are nonisolated because for any 
$\alpha \in \mathcal{I}$ and any $\epsilon>0$ 
we can find a geodesic of infinite type whose cutting sequence diverges
from the one of $\gamma_\alpha$ only after the $N_\epsilon$th term.
On the other hand, $\mathcal{I}$ is discrete, because within the
$\epsilon$-neighbourhood of $\gamma_\alpha$ we can construct a geodesic
which fills $D_n$ in finite time: just choose it to have a cutting 
sequence which agrees with the one of $\gamma_\alpha$ for $N_\epsilon$ terms,
and to then career off along some path which decomposes $D_n$ into
disks and once-punctured disks.   

Finally, the last part of theorem 3.3(c) holds because each of the 
countably many elements of $B_n$ fixes only a countable number of points
$\alpha \in (0,\pi)$ with the property that $\gamma_\alpha$ fills $D_n$. 
In order to see this, we note that
for {\it irreducible} elements of $B_n$ theorem 5.5 of \cite{CB} states that
there is only a finite number of fixed points on the circle at infinity. 
If an element $\phi$ of $B_n$ is {\it reducible}, then we leave it to the
reader to check that the result follows from the following facts:
(1) One can find a maximal invariant system $C$ of disjoint embedded arcs and 
circles in $D_n$.
\ (2) If $\phi$ acts nontrivially
on a component of $D_n \basl C$ which is cut in a nontrivial way by a 
{\it finite} segment of $\gamma_\alpha$, 
then it acts nontrivially on $\gamma_\alpha$
(for if it didn't then the collection $C$ would not be maximal). 
\ (3) A geodesic $\gamma_\alpha$ that fills $D_n$ has to enter every 
component of $D_n\basl C$ at least once, and $\phi$ acts nontrivially 
either on the first or, failing that, on the second component of 
$\gamma_\alpha\cap (D_n \basl C)$ (because it cannot act trivially on two
adjacent components of $D_n \basl C$). 
\ (4) There is a countable infinity of isotopy classes of embedded arcs 
from the basepoint of $D_n$ to $C$. 
\end{proof}

We recall from the beginning of the section that to every geodesic 
$\gamma_\alpha$ of infinite type we have associated a ``critical disk''
$D_c$ which contains most of the last arc of $C(\gamma_\alpha)$.
The fundamental property of geodesics of infinite type which we shall
use several times is the following.
\begin{lem}\sl For any geodesic of infinite type $\gamma_\alpha$ and
for any $\epsilon>0$ there exists a geodesic $\gamma_{\alpha^+}$  with
$\alpha^+ \in (\alpha, \alpha+\epsilon)$ such that $\gamma_{\alpha^+}$
falls into a puncture and has no self-intersections inside $D_c$.
\end{lem}
\begin{proof} It suffices to prove the lemma in the special case $D_c=D_n$, 
ie when the geodesic $\gamma_\alpha$ is embedded.  
We suppose, for a contradiction, that there exists an
$\epsilon>0$ such that no $\gamma_\beta$ with 
$\beta \in (\alpha, \alpha+\epsilon)$ is embedded {\it and} falls into
a puncture. Our aim is to reach the contradiction that $\gamma_\alpha$
ends in an infinite spiral.

We continue to use the notions concerning cutting sequences
introduced above: we choose arbitrarily a fundamental domain, and
we shall denote by $\gamma_\alpha^k$ the initial segment
of $\gamma_\alpha$ up to its $k$th intersection with the boundary of
the fundamental domain.
We recall that, given $\gamma_\alpha$ and $\epsilon>0$, 
we can find an $N=N_\epsilon\in \N$ such that any
geodesic $\gamma_\beta$ with $\gamma_\beta^N=\gamma_\alpha^N$
satisfies $|\alpha - \beta|<\epsilon$.
We now consider the arc $\gamma_\alpha^{N+1}$: it ends on some boundary
arc of the fundamental domain which we denote $a$.
The orientation of $\gamma_\alpha$ gives rise to
a notion of the part of $a$ ``to the left'' and ``to the right of'' the
end point of $\gamma_\alpha^{N+1}$. The arc $\gamma_\alpha^{N+1}$  has an
intersection with the interior of the ``left'' part of $a$, for if hadn't 
we could obtain an embedded arc $\gamma_\beta$ with 
$\beta \in (\alpha, \alpha+\epsilon)$ by adjoining to the end point of 
$\gamma_\alpha^N$  an arc falling into the puncture at the left end of $a$;
this would contradict the hypothesis.
Thus it makes sense to define $\Gamma \subseteq D_n$
to be the union of $\gamma_\alpha^{N+1}$ and a segment of $a$ from the 
end point of 
$\gamma_\alpha^{N+1}$ to the left, up to the next intersection with
$\gamma_\alpha^{N+1}$ (see figure \ref{gammasp}).

\begin{figure}[htb] %
\cl{
\relabelbox
\epsfbox{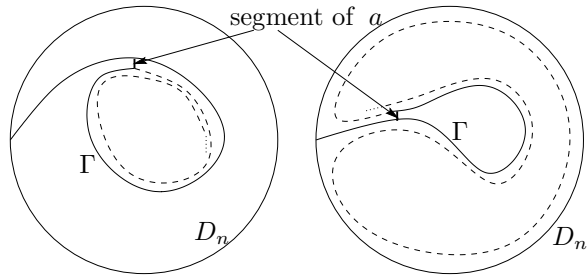} 
\relabel{Gamma}{\small $\Gamma$}
\relabel{Gamma2}{\small $\Gamma$}
\relabel{Dn}{\small $D_n$}
\relabel{Dn2}{\small $D_n$}
\relabel{segment of a}{\small segment of \thinspace$a$}
\endrelabelbox
}
\caption{The two possible shapes of $\Gamma$, and (dashed) the
resulting geodesic $\gamma_\alpha$}
\label{gammasp}
\end{figure}

We now observe that $D_n \basl \Gamma$ has two path components, each
containing at least one puncture; moreover, $\gamma_\alpha$ cannot
intersect any geodesic arc connecting two punctures in the same
component, because the first time it did we could drop it into the
puncture at the left end of the arc and obtain a contradiction as
before. It follows that $\gamma_\alpha$ has to spiral along the boundary
of one of the components of $D_n \basl \Gamma$.
\end{proof}

\begin{prop}\sl All orderings, even partial ones, arising from geodesics
$\gamma_\alpha$ of infinite type are non-discrete.
\end{prop}
\begin{proof} We shall prove the following stronger statement: 
for any $\epsilon>0$ there exists an element $\phi \in \MCG(D_n)=B_n$ 
such that $\phi(\alpha)\in (\alpha,\alpha+\epsilon)$.

We choose $\alpha^+$ 
as in the previous lemma. We consider the boundary curve $\tau$ of
a regular neighbourhood of $\partial D_c \cup \gamma_{\alpha^+}$ in
$D_c$. This curve $\tau$ is disjoint from $\gamma_{\alpha^+}$, while
any curve isotopic to $\tau$ necessarily intersects $\gamma_\alpha$.
Thus for the positive Dehn twist $T$ along $\tau$ we have that 
$T(\alpha)>\alpha$ (by Proposition 2.4), and that $T(\alpha^+)=\alpha^+$.
It follows that $T(\alpha) \in (\alpha, \alpha^+) \subseteq 
(\alpha, \alpha+\epsilon)$.
\end{proof}

\begin{proof}[Proof of Theorem 3.4(a)] 
Given a geodesic $\gamma_\alpha$ of finite, and a geodesic $\gamma_\beta$ 
of infinite type, our aim is to prove that $\gamma_\alpha$ and 
$\gamma_\beta$ cannot induce the same orderings of $B_n$. 

As seen in corollary 6.2, orderings arising from geodesics which fill 
the surface in finite time are the same as orderings arising from total 
curve diagrams, which are discrete by lemma 4.5. By contrast, we have 
from proposition 7.2 that infinite type orderings are not discrete. This 
proves the theorem in the special case where the finite type geodesic
fills the surface.

In the case where the finite type geodesic $\gamma_\alpha$ does {\it not} 
fill the surface, we consider the subsurface 
$D_\alpha := D_n \basl NC(\gamma_\alpha)$, ie the maximal subsurface with 
geodesic boundary which is disjoint from $\gamma_\alpha$. We observe that 
$D_\alpha$ is a disjoint union of disks, each containing at least two 
punctures. Any homeomorphism $\phi$ of $D_n$ with support in $D_\alpha$ has 
the property that $\phi(\alpha)=\alpha$.
 
If $D_\alpha \cap \gamma_\beta \neq \emptyset$ then there exists a 
homeomorphism $\phi$ with support in $D_\alpha$ such that 
$\phi(\beta)\neq \beta$, and it follows that the orderings induced by 
$\alpha$ and $\beta$ are different. 

If, on the other hand, $D_\alpha \cap \gamma_\beta = \emptyset$,
then we squash each component of $D_\alpha$ to a puncture; the result is
a disk with say $m$ punctures, where $m<n$, which we denote $D_m$.
We now consider the subgroup $B_m^P$ of $B_m = \MCG(D_m)$ of all mapping 
classes which fix those punctures of $D_m$ that came from squashed 
components of $D_\alpha$. This is a finite index subgroup of $B_m$, and
the orderings of $B_n$ determined by $\alpha$ and $\beta$ induce quotient 
orderings on $B_m^P$.
Another way to describe these quotient orderings is to repeat the 
Thurston-construction for the disk $D_m$: one can equip $D_m$ with a 
hyperbolic metric, and then the geodesics $\gamma_\alpha$ and $\gamma_\beta$ 
project to quasigeodesics in $D_m$. These quasigeodesics determine points 
at infinity of the universal cover of $D_m$, and hence give rise to
orderings of $B_m$.

The geodesic in $D_m$ which is homotopic to the projection of $\gamma_\alpha$
is again of finite type; the crucial observation now is that it fills
$D_m$, so that the quotient ordering on $B_m^P$ is discrete by Lemma 4.5.
Similarly, a geodesic in $D_m$ homotopic to the projection of $\gamma_\beta$
is again of infinite type, and hence induces a non-discrete ordering, by
Proposition 7.2. So the $\alpha$- and $\beta$-orderings on $B_n$ give 
rise to different quotient orderings on $B_m^P$, and are therefore different.
\end{proof}

As seen above, every geodesic of infinite type gives rise to a curve diagram
``of infinite type'', which is like a curve diagram of finite type, except
that the arc with maximal label is, up to isotopy, an infinite geodesic 
which does not fall into a puncture or a spiral. All but a finite 
initial segment of this arc lies in the ``critical disk'' $D_c$.
There is an obvious generalisation of the notion of loose isotopy:

\begin{definition}\sl  Two curve diagrams of infinite type are 
{\it loosely isotopic} if they are related by 
(1) continuous deformation, ie a path in the space of all curve diagrams 
of infinite type; and (2) pulling loops around punctures tight. 
\end{definition}

This is exactly the same as in the finite type case, except that no
``pulling loops around punctures tight''-procedure is defined for the
last arc. We are now ready to state and prove the main classification 
theorem for orderings of $B_n$ of infinite type.
\begin{theorem}
Two geodesics $\gamma_\alpha$ and $\gamma_\beta$ of infinite type give
rise to the same ordering of $B_n$ if and only if their associated
curve diagrams $C(\gamma_\alpha)$ and $C(\gamma_\beta)$ are loosely
isotopic.
\end{theorem}
\begin{proof} By the results in the previous sections, it suffices to 
prove that two {\it embedded} geodesics $\gamma_\alpha$ and $\gamma_\beta$ 
of infinite type give rise to the same ordering of $B_n$ if and only if
$\beta=\Delta^{2k}(\alpha)$ for some $k\in \Z$, ie if $\gamma_\alpha$ and 
$\gamma_\beta$ are related by a slide of the starting point around
$\partial D_n$. 

The implication ``$\Leftarrow$'' is clear. Conversely, for the implication
``$\Rightarrow$'', we suppose that $\gamma_\alpha$ and $\gamma_\beta$ are 
not related by a slide of the starting point, and without loss of generality
we say $\alpha > \beta$. Our aim is to construct a homeomorphism which
is positive in the $\alpha$- and negative in the $\beta$-ordering, ie
which sends $\alpha$ ``more to the left'' and $\beta$ ``more to the right''.
Our argument will be a refinement of the proof of the implication 
``$\Rightarrow$'' of 5.2(a). 

By lemma 7.1 we can construct embedded geodesics $\gamma_{\alpha^+}$
and $\gamma_{\beta^+}$ which fall into punctures, and lie an arbitrarily 
small amount to the left of $\gamma_\alpha$ respectively $\gamma_\beta$.
We define the curves $\tau_{\alpha^+}$ and 
$\tau_{\beta^+}$ to be the geodesic representatives of the boundary curves 
of regular neighbourhoods in $D_n$ of $\partial D_n \cup \gamma_{\alpha^+}$ 
and ${\partial D_n \cup \gamma_{\beta^+}}$ respectively. We denote by
$T_{\alpha^+}$ respectively $T_{\beta^+}$ the positive Dehn twists
along these curves. Our desired homeomorphism will be of the form 
$T_{\alpha^+}^{-k} \circ T_{\beta^+}$, with carefully chosen values
of $\alpha^+$ and $\beta^+$, and $k\in \N$ very large.

We also define the two-sided infinite geodesic $\tau_\alpha$ to be the
geodesic which is disjoint from $\gamma_\alpha$, and isotopic to the
boundary of a neighbourhood of $\gamma_\alpha \cup \partial D_n$ in $D_n$.
More formally, in the universal cover $D_n^\eqsim$ we consider two
liftings of $\gamma_\alpha$, namely $\widetilde{\gamma}_\alpha$ (which 
starts at the basepoint of $D_n^\eqsim$), and the lifting whose starting
point also lies on $\Pi$, and is obtained from the basepoint of 
$\widetilde{D}_n$ by
lifting the path in $D_n$ once around $\partial D_n$. The end points of
these geodesics lie on the circle at infinity, and $\tau_\alpha$ is just
the projection of the geodesic connecting them.

Since $\gamma_\alpha$ and $\gamma_\beta$ are not loosely isotopic, we
have that $\gamma_\beta$ intersects $\tau_\alpha$. By choosing 
$\beta^+$ sufficiently close to $\beta$ we can now achieve that the
initial segments of $\gamma_\beta$ and $\gamma_{\beta^+}$ up to their 
first point of intersection with $\tau_\alpha$ are isotopic with end points 
sliding in $\tau_\alpha$. This gives our choice of $\beta^+$, and it remains
to choose $\alpha^+$ and $k$.

The crucial observation concerning $\tau_\alpha$ is that it can be
arbitrarily closely approximated by the curves $\tau_{\alpha^+}$, by 
choosing $\alpha^+$ sufficiently close to $\alpha$. More precisely, 
in the universal cover $D_n^\eqsim$ we consider the preimages 
of $\tau_\alpha$ and of $\tau_{\alpha^+}$. Each of them has
infinitely many path components; we choose one distinguished component
for each, namely the first ones that $\gamma_\beta$ intersects.
Our observation now is that as $\alpha^+$ tends to $\alpha$, 
the end points of the distinguished component of the preimage of
$\tau_{\alpha^+}$ tend to the end points of the distinguished component
of the preimage of $\tau_\alpha$.

We now turn to the choice of $\alpha^+$. By proposition 2.4 we have that 
$T_{\beta^+}(\alpha)>\alpha$. 
By lemma 7.1 we can now choose $\alpha^+$ close to $\alpha$ such that 
$T_{\beta^+}(\alpha)>\alpha^+>\alpha$. By eventually pushing $\alpha^+$
even closer to $\alpha$, we can in addition insist (by the observation
concerning $\tau_\alpha$ above), that the initial 
segments of $\gamma_\beta$ and $\gamma_{\beta^+}$ 
up to their first point of intersection with $\tau_{\alpha^+}$ are also 
isotopic with end points sliding in $\tau_{\alpha^+}$. This gives our choice
of $\alpha^+$. 

We have arrived at the following setup: we have the three points
$\beta^+=T_{\beta^+}(\beta^+)>T_{\beta^+}(\beta)>\beta$ in 
$\partial D_n^\eqsim \basl \Pi$, and they all lie between
the two end points $\delta_l$ and $\delta_r$ of the distinguished 
lifting of $\tau_{\alpha^+}$ (here the indices $l$ and $r$ stand for
``left'' and ``right'', so $\delta_l>\delta_r$). For any point $\delta$
with $\delta_l>\delta>\delta_r$ we consider the action of the positive
Dehn twist $T_{\alpha^+}$ on the geodesic $\gamma_\delta$. We observe
that the limit $\lim_{k\to \infty} T_{\alpha^+}^{-k}(\delta) = \delta_r$.
In particular for $\delta := \beta^+$ it follows that for sufficiently 
large $k$ we have $T_{\alpha^+}^{-k}(\beta^+)<\beta$.
This gives our choice of $k$.

To summarise, we have 
$$ T_{\alpha^+}^{-k} \circ T_{\beta^+}(\alpha)>T_{\alpha^+}^{-k}(\alpha^+)
=\alpha^+>\alpha$$ and 
$$ T_{\alpha^+}^{-k} \circ T_{\beta^+}(\beta)<
 T_{\alpha^+}^{-k} \circ T_{\beta^+}(\beta^+)=
 T_{\alpha^+}^{-k}(\beta^+)<\beta,
$$
ie $T_{\alpha^+}^{-k} \circ T_{\beta^+}$ is positive in
the $\alpha$-, but negative in the $\beta$-ordering.
\end{proof}
\begin{proof}[Proof of Theorem 3.4(c)] This is an immediate consequence
of Theorem 7.4 \end{proof}

{\bf Acknowledgements } We are very grateful to W.P.~Thurston for sharing 
with us his idea for the construction of orderings on $\MCG(S)$. We also
thank the following people for helpful comments: S.~Chmutov, D.~Cooper, 
L.~Paris, C.~Rourke, L.~Rudolph, and M.~Scharlemann.


\end{document}